%% file: motion_di.tex
\documentclass[journal]{IEEEtran}

\usepackage{verbatim}
\usepackage{cite}
\usepackage[pdftex]{graphicx}
\graphicspath{{./images/}}

\usepackage[cmex10]{amsmath}
\usepackage{amssymb}
\usepackage{latexsym}
\usepackage{psfrag}
\usepackage{epsfig}
\usepackage{graphicx}
\usepackage{epstopdf}
\usepackage{multicol}
\usepackage{arydshln}
\usepackage{color}
\usepackage{tikz}
\usetikzlibrary{shapes,arrows}
\usetikzlibrary{calc}

\newtheorem{theorem}{\textbf{Theorem}}[section]
\newtheorem{lemma}[theorem]{\textbf{Lemma}}
\newtheorem{proposition}[theorem]{\textbf{Proposition}}

\newtheorem{remark}[theorem]{Remark}

\newtheorem{assumption}[theorem]{\textbf{Assumption}}

\newcommand{\R}{{\rm I\!R}}
\newcommand{\dfb}{\stackrel{\Delta}{=}}

\usepackage{caption}
\usepackage{subcaption}

\usepackage{url}

\tikzstyle{block} = [draw, fill=white, rectangle, 
	minimum height=3em, minimum width=6em, text width=5em, text centered]
\tikzstyle{sum} = [draw, fill=white, circle, node distance=2cm]
\tikzstyle{input} = [coordinate]
\tikzstyle{output} = [coordinate]
\tikzstyle{pinstyle} = [pin edge={to-,thin,black}]
\tikzstyle{branch} = [circle,inner sep=0pt,minimum size=1mm,fill=black,draw=black]

\tikzstyle{vertex}=[circle,fill=black!25,minimum size=20pt,inner sep=0pt]
\tikzstyle{selected vertex} = [vertex, fill=red!24]
\tikzstyle{edge} = [draw,thick,-]
\tikzstyle{dedge} = [draw,thick,->]
\tikzstyle{shadowdedge} = [draw, dotted,->]
\tikzstyle{weight} = [font=\small]
\tikzstyle{selected edge} = [draw,line width=5pt,-,red!50]
\tikzstyle{ignored edge} = [draw,line width=5pt,-,black!20]

\begin{document}
\title{Taming mismatches in inter-agent distances for the formation-motion control of second-order agents.}
%
%
\author{Hector~Garcia de Marina,~\IEEEmembership{Member,~IEEE,}
	Bayu~Jayawardhana,~\IEEEmembership{Senior Member,~IEEE,}
        and Ming~Cao,~\IEEEmembership{Senior Member,~IEEE}%
\thanks{H. Garcia de Marina is with the Ecole Nationale de l'Aviation Civile, University of Toulouse, Toulouse, France. (e-mail: hgdemarina@ieee.org).}
\thanks{B. Jayawardhana and M. Cao are with the Engineering and Technology Institute of Groningen, University of Groningen, 9747 AG Groningen, The Netherlands. (e-mail: \{b.jayawardhana, m.cao\}@rug.nl). This work was supported by the the EU INTERREG program under the auspices of the SMARTBOT project and the work of Cao was also supported by the European Research Council (ERC-StG-307207).}}


\markboth{IEEE TRANSACTIONS ON AUTOMATIC CONTROL}%
{}

\maketitle

\begin{abstract}
This paper presents the analysis on the influence of distance mismatches on the standard gradient-based rigid formation control for second-order agents. It is shown that, similar to the first-order case as recently discussed in the literature, these mismatches introduce two undesired group behaviors: a distorted final shape and a steady-state motion of the group formation. We show that such undesired behaviors can be eliminated by combining the standard formation control law with distributed estimators. Finally, we show how the mismatches can be effectively employed as design parameters in order to control a combined translational and rotational motion of the formation.
\end{abstract}


\section{Introduction}
\IEEEPARstart{R}{ecent} years have witnessed a growing interest in coordinated robot tasks, such as, area exploration and surveillance \cite{7139494, 844100}, robot formation movement for energy efficiency \cite{tsugawa2011automated}, and tracking and enclosing a target \cite{guo2010local,hara2008distributed}. In many of these team-work scenarios, one of the key tasks for the agents is to form and maintain a prescribed formation shape. Gradient-based control has been widely used for this purpose \cite{BaArWe11}. In particular, \emph{distance-based} control \cite{OlMu02,KrBrFr08,YuAnDaFi09,CaYuAn11,oh2015survey} has gained popularity since the agents can work with their own local coordinates and the desired shape of the formation under control is \emph{exponentially stable} \cite{SuLiAn15,SuAn15}. However, exponential stability cannot prevent undesired steady-state collective motions, e.g., constant drift, if disturbances, such as biases in the range sensors or equivalently mismatches between the prescribed distances for neighboring agents, are present. This misbehavior has been carefully studied for agents governed by single integrator dynamics in \cite{MouMorseBelSunAnd15,Hem14} and an effective solution to get rid of such misbehavior using estimators has been reported in \cite{MarCaoJa15}. The proposed estimator-based tool works without requiring any communication among agents. This is a desired feature since the above mentioned issue with the mismatches cannot usually be solved in a straightforward way by sending back and forth more communication information between the agents for several reasons: the sensing bias can be time-varying due to factors such as the environment's temperature; the same range sensor can produce different readings for the same physical distance in face of random measurement noises, making communication-based correction costly; a continuous (or regular) communication among agents may not even be possible or desired; and the agents may have different clocks complicating the data comparison in real-time computation.

In this paper for second-order agents, i.e., agents modeled by double integrator dynamics, we extend the recent findings in \cite{MouMorseBelSunAnd15,Hem14, MarCaoJa15, MaJaCa15} on \emph{mismatched} formation control, the cancellation of the undesired effects via distributed estimation, and the formation motion control by turning the mismatches into distributed parameters. There are practical benefits justifying such extension. For example, a formation controller employing undirected sensing topologies has inherent stability properties that are not present if \emph{directed} topologies\footnote{For directed topologies, only one agent per edge controls the inter-agent distance. Hence, it is free of mismatches.} are employed, especially for agents with higher order dynamics \cite{BaArWe11}. Another key aspect of employing second-order dynamics is that the control actions can be used as the \emph{desired acceleration} in a guidance system feeding the tracking controller of a mechanical system, such as the ones proposed for quadrotors in \cite{MeNaVi12,smeur2016gust} or marine vessels in \cite{fossen2002marine}. This clearly simplifies such tracking controllers compared to the case of only providing \emph{desired velocities} derived from first-order agent dynamics. 

The extensions presented in this paper require new technical constructions that go much beyond what is needed for first-order agents. For example, a key step in the logic presented in \cite{MouMorseBelSunAnd15,Hem14,MaJaCa15,MarCaoJa15} is based on the fact that the dynamics of the error signal, which measures the distortion with respect to the desired shape, is an autonomous system for first-order agents. This does not hold anymore in the second-order case. Consequently, additional technical steps are developed for extending the results of a \emph{mismatched} formation control to second-order agents.

The problem of motion and formation control can be solved simultaneously if the mismatches are well understood and not treated as disturbances but as design parameters. Indeed, this is the strategy followed in \cite{MaJaCa15}, where we have turned mismatches in the prescribed distances into design parameters in order to manipulate the way how the collective motion is realized. The act of turning mismatches to distributed motion parameters allows one to address more complicated problems such as moving a rigid formation (including rotational motion) without leaders \cite{BaArWe11}, or tracking and enclosing a free target. The desired motion is designed with respect to a frame of coordinates fixed at the desired \emph{rigid-body} shape. This design allows to preserve the ordering of the agents during the motion, e.g., which agent is leading at the front of the group. An important aspect of this approach is that for \emph{low} motion speeds the agents do not need to measure any relative velocities as what is commonly required in swarms of second-order agents \cite{deghatcombined}. Consequently, the sensing requirements are reduced. It will be shown that the motion parameters for second-order agents can be directly computed from the first-order case as in \cite{MaJaCa15}. Unfortunately, the Lyapunov function for proving the stability of the desired steady-state motion is not as straightforward as in \cite{MaJaCa15}, since the \emph{standard} quadratic function, involving the norms of the errors regarding the distortion of the desired shape and the velocity of the agents, can only be used to prove asymptotic \cite{BaArWe11} but not exponential stability, which is the key property for determining the necessary gains and region of attraction for the formation and motion controller in \cite{MaJaCa15}. Nevertheless, we will be able to show that the stability of the closed-loop formation motion system is indeed exponential for second-order agents too.

The rest of the paper is organized as follows. After some definitions and clarification of notation in Section \ref{sec: pre}, we show in Section \ref{sec: issues} the analysis about the undesired effects on a formation of second-order agents in the presence of small mismatches. It turns out that in addition to a distortion with respect to the desired shape, we have an undesired steady-state collective motion which is consistent with the one described in \cite{SunMouAndMor14} for first-order agents. In Section \ref{sec: esti}, we propose two designs for a distributed estimator in order to prevent the mentioned robustness issues. The first design eliminates the undesired steady-state motion and is able to bound the steady-state distortion with respect to a generic rigid desired shape. In particular, this first design completely removes the distortion for the special cases of triangles and tetrahedrons. The second design is less straightforward and requires the calculations of some lower bounds for certain gains. However, it eliminates the undesired steady-state motion and distortion for any desired shape. In Section \ref{sec: motion}, we turn the mismatches into distributed motion parameters in a similar way as in \cite{MaJaCa15} in order to design the stationary motion of the formation without distorting the desired shape. Finally, in order to validate the proposed algorithms from Sections \ref{sec: esti} and \ref{sec: motion}, simulation results with second-order agents are provided in Section \ref{sec: exp}.

\section{Preliminaries}
\label{sec: pre}
In this section, we introduce some notations and concepts related to graphs and rigid formations. For a given matrix $A\in\R^{n\times p}$, define $\overline A \dfb A \otimes I_m \in\R^{nm\times pm}$, where the symbol $\otimes$ denotes the Kronecker product, $m = 2$ for the 2D formation case or $m=3$ for the 3D one, and $I_m$ is the $m$-dimensional identity matrix. For a stacked vector $x\dfb \begin{bmatrix}x_1^T & x_2^T & \dots & x_k^T\end{bmatrix}^T$ with $x_i\in\R^{n}, i\in\{1,\dots,k\}$, we define the diagonal matrix $D_x \dfb \operatorname{diag}\{x_i\}_{i\in\{1,\dots,k\}} \in\R^{kn\times k}$. We denote by $|\mathcal{X}|$ the cardinality of the set $\mathcal{X}$ and by $||x||$ the Euclidean norm of a vector $x$. We use $\mathbf{1}_{n\times m}$ and $\mathbf{0}_{n\times m}$ to denote the all-one and all-zero matrix in $\R^{n\times m}$ respectively and we will omit the subscript if the dimensions are clear from the context.

\subsection{Graphs and Minimally Rigid Formations}
\label{sec: preA}
We consider a formation of $n\geq 2$ autonomous agents whose positions are denoted by $p_i\in\R^m$. The agents can measure their relative positions with respect to its neighbors. This sensing topology is given by an undirected graph $\mathbb{G} = (\mathcal{V}, \mathcal{E})$ with the vertex set $\mathcal{V} = \{1, \dots, n\}$ and the ordered edge set $\mathcal{E}\subseteq\mathcal{V}\times\mathcal{V}$. The set $\mathcal{N}_i$ of the neighbors of agent $i$ is defined by $\mathcal{N}_i\dfb\{j\in\mathcal{V}:(i,j)\in\mathcal{E}\}$. We define the elements of the incidence matrix $B\in\R^{|\mathcal{V}|\times|\mathcal{E}|}$ for $\mathbb{G}$ by
\begin{equation*}
	b_{ik} \dfb \begin{cases}+1 \quad \text{if} \quad i = {\mathcal{E}_k^{\text{tail}}} \\
		-1 \quad \text{if} \quad i = {\mathcal{E}_k^{\text{head}}} \\
		0 \quad \text{otherwise,}
	\end{cases}
\end{equation*}
where $\mathcal{E}_k^{\text{tail}}$ and $\mathcal{E}_k^{\text{head}}$ denote the tail and head nodes, respectively, of the edge $\mathcal{E}_k$, i.e. $\mathcal{E}_k = (\mathcal{E}_k^{\text{tail}},\mathcal{E}_k^{\text{head}})$. A \emph{framework} is defined by the pair $(\mathbb{G}, p)$, where $p = \operatorname{col}\{p_1, \dots, p_n\}$ is the stacked vector of the agents' positions $p_i,i\in\{1,\dots,n\}$. With this last definition at hand, we define the stacked vector of the measured relative positions by
\begin{equation*}
	z = \overline B^Tp,
\end{equation*}
where each vector $z_k = p_i - p_j$ in $z$ corresponds to the relative position associated with the edge $\mathcal{E}_k = (i, j)$.

For a given stacked vector of desired relative positions $z^* = [\begin{smallmatrix}{z_1^*}^T & {z_2^*}^T & \dots & {z_{|\mathcal{E}|}^*}^T\end{smallmatrix}]^T$, the resulting set $\mathcal{Z}$ of the possible formations with the same shape is defined by
\begin{equation}
	\mathcal{Z} \dfb \left \{\left(I_{|\mathcal{E}|} \otimes \mathcal{R}\right)z^* \right \}, \label{eq: Zp}
\end{equation}
where $\mathcal{R}$ is the set of all rotational matrices in 2D or 3D. Roughly speaking, $\mathcal{Z}$ consists of all those formation positions that are obtained by rotating $z^*$.

\begin{figure}
	\centering
	\begin{subfigure}{0.2\columnwidth}
		\input{FIGsquaredef.tex}
		\caption{}
	\end{subfigure} \quad\quad
	\begin{subfigure}{0.2\columnwidth}
		\input{FIGsquareRig.tex}
		\caption{}
	\end{subfigure} \quad 
	\begin{subfigure}{0.2\columnwidth}
		\input{FIGtetra.tex}
		\caption{}
	\end{subfigure}
	\caption{a) The square without an inner diagonal is not rigid since we can smoothly move the top two nodes while keeping the other two fixed without breaking the distance constraints; b) The square can be done locally minimally rigid in $\R^2$ if we add an inner diagonal; c) The tetrahedron in $\R^3$ is infinitesimally and minimally rigid.}
	\label{fig: rig}
\end{figure}
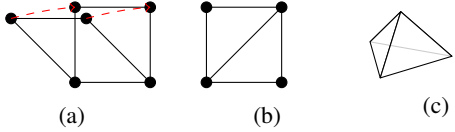

Let us now briefly recall the notions of infinitesimally rigid framework and minimally rigid framework from \cite{AnYuFiHe08}. Define the edge function $f_\mathbb{G}$ by $f_{\mathbb{G}}(p) = \mathop{\text{col}}\limits_{k}\big(\|z_k\|^2\big)$ and we denote its Jacobian by
\begin{equation}
2R(z) = 2D_z^T\overline B^T,
\end{equation}
where $R(z)$ is called the {\it rigidity matrix}. A framework $(\mathbb{G}, p)$ is {\it infinitesimally rigid} if $\text{rank}\{R(z)\} = 2n-3$ when it is embedded in $\mathbb{R}^2$ or if $\text{rank}\{R(z)\} = 3n-6$ when it is embedded in $\mathbb{R}^3$. Additionally, if $|\mathcal E|=2n-3$ in the 2D case or $|\mathcal E|=3n-6$ in the 3D case, then the framework is called {\it minimally rigid}. Roughly speaking, the only motions that one can perform over the agents in an infinitesimally and minimally rigid framework, while they are already in the desired shape, are the ones defining translations and rotations of the whole shape. Some graphical examples of infinitesimally and minimally rigid frameworks are shown in Figure \ref{fig: rig}.
If $(\mathbb G, p)$ is infinitesimally and minimally rigid, then, similar to the above, we can define the set of resulting formations $\mathcal{D}$ by
\begin{align}
\mathcal{D} \dfb & \Big \{z  \, | \, ||z_k||=d_k , k\in\{1, \dots, |\mathcal{E}|\} \Big \}, \nonumber
\end{align}
where $d_k = ||z_k^*||, k\in\{1, \dots, |\mathcal{E}|\}$.

Note that in general it holds that $\mathcal{Z}\subseteq\mathcal{D}$. For a desired shape, one can always design $\mathbb G$ to make the formation infinitesimally and minimally rigid. In fact, an infinitesimally and minimally rigid framework with two or more vertices in $\R^2$ can always be constructed through the Henneberg construction \cite{Hen11}. In $\R^3$ one can construct a set of infinitesimally and minimally rigid frameworks via insertion starting from a tetrahedron, if each newly added vertex with three newly links forms another tetrahedron as well.
\subsection{Frames of coordinates}
It will be useful for describing the motions of the infinitesimally and minimally rigid formation to define a frame of coordinates fixed to the desired formation itself. We denote by $O_g$  the \emph{global frame} of coordinates fixed at some point of $\R^m$ with some arbitrary fixed orientation. In a similar way, we denote by $O_b$  the \emph{body frame} fixed at the centroid $p_c$ of the desired rigid formation. Furthermore, if we rotate the rigid formation with respect to $O_g$, then $O_b$ is also rotated in the same manner. Let $^bp_j$ denote the position of agent $j$ with respect to $O_b$. To simplify the notation whenever we represent an agents' variable with respect to $O_g$, the superscript is omitted, i.e., $p_j \dfb {^gp_j}$.

\section{Robustness issues due to mismatches in formation gradient-based control}
\label{sec: issues}
\subsection{Gradient Control}
Consider a formation of $n$ agents with the sensing topology $\mathbb{G}$ for measuring the relative positions among the agents. The agents are modelled by a second-order system given by
\begin{equation}
	\begin{cases}
	\dot p = v \\
	\dot v = u,
	\end{cases}
	\label{eq: pdyn}
\end{equation}
where $u$ and $v$ are the stacked vector of control inputs $u_i\in\R^m$ and vector of agents' velocity $v_i\in\R^m$ for $i=\{1,\dots,n\}$ respectively.

In order to control the shape, for each edge $\mathcal{E}_k = (i,j)$ in the infinitesimally and minimally rigid framework we assign the following potential function $V_k$
\begin{equation*}
	V_k(||z_k||) = \frac{1}{4}(||z_k||^2 - d_k^2)^2,
\end{equation*}
with the gradient along $p_i$ or $p_j$ given by 
\begin{equation*}
	\nabla_{p_i}V_k = -\nabla_{p_j}V_k = z_k (||z_k||^2 - d_k^2).
\end{equation*}
In order to control the agents' velocities, for each agent $i$ in the infinitesimally and minimally rigid framework we assign the following (kinetic) energy function $T_i$
\begin{equation*}
	T_i(v_i) = \frac{1}{2}||v_i||^2,
\end{equation*}
with the gradient along $v_i$ be given by
\begin{equation*}
	\nabla_{v_i}T_i = v_i.
\end{equation*}
One can check that for the potential function
\begin{equation}
	\phi(p, v) = \sum_{i=1}^{|\mathcal{V}|} T_i + \sum_{k=1}^{|\mathcal{E}|}V_k,
	\label{eq: phi}
\end{equation}
the closed-loop system (\ref{eq: pdyn}) with the control input
\begin{equation}
	u = -\nabla_v \phi -\nabla_p \phi, \label{eq: un}
\end{equation}
becomes the following dissipative Hamiltonian system \cite{schaft2006port}
\begin{equation}
	\begin{cases}
	\dot p = \nabla_v \phi \\
	\dot v = -\nabla_v \phi -\nabla_p \phi.
	\end{cases}
	\label{eq: H}
\end{equation}
Considering (\ref{eq: phi}) as the storage energy function of the Hamiltonian system (\ref{eq: H}), one can show the local asymptotic convergence of the formation to the shape given by $\mathcal{D}$ and all the agents' velocities to zero \cite{BaArWe11,Oh14}.

Consider the following one-parameter family of dynamical systems $\mathcal{H}_\lambda$ given by
\begin{equation}
\begin{bmatrix}\dot p \\ \dot v\end{bmatrix} = -
\begin{bmatrix}\lambda I_{m|\mathcal{V}|} & -(1-\lambda)I_{m|\mathcal{V}|} \\
	(1-\lambda)I_{m|\mathcal{V}|} & I_{m|\mathcal{V}|}
\end{bmatrix}\begin{bmatrix}\nabla_p \phi \\ \nabla_v \phi \end{bmatrix},
	\label{eq: Hl}
\end{equation}
where $\lambda \in [0, 1]$, which defines all convex combinations of the Hamiltonian system (\ref{eq: H}) and a gradient system. The family $\mathcal{H}_\lambda$ has two important properties summarized in the following lemma.
\vspace{5pt}
\begin{lemma}\cite{Oh14}
	\label{lem: H}
\begin{itemize}
\item For all $\lambda \in [0, 1]$, the equilibrium set of $\mathcal{H}_\lambda$ is given by the set of the critical points of the potential function $\phi$, i.e. $E_{p,v} = \left\{\begin{bmatrix}p^T & v^T\end{bmatrix}^T : \nabla\phi = \mathbf{0}\right\}$.
	\item For any equilibrium $\begin{bmatrix}p^T & v^T\end{bmatrix}^T\in E_{p,v}$ and for all $\lambda \in [0, 1]$, the numbers of the stable, neutral, and unstable eigenvalues of the Jacobian of $\mathcal{H}_\lambda$ are the same  and independent of $\lambda$.
\end{itemize}
\end{lemma}

This result has been exploited in \cite{SuAn15} in order to show the local exponential convergence of $z(t)$ and $v(t)$ to $\mathcal{D}$ and $\mathbf{0}$ respectively. In the following brief exposition we revisit such exponential stability via a combination of Lyapunov argument and Lemma \ref{lem: H}, which will play an important role in Section \ref{sec: is2}.

Define the distance error corresponding to the edge $\mathcal{E}_k$ by
\begin{equation}
e_k = ||z_k||^2 - d^2_k,
	\label{eq: error2}
\end{equation}
whose time derivative is given by $\dot e_k = 2z_k^T\dot{z_k}$. Consider the following autonomous system derived from (\ref{eq: Hl}) with
$\lambda = 0.5$
\begin{align}
	\dot p &= -\frac{1}{2}\overline BD_ze + \frac{1}{2}v \nonumber \\
	\dot z &= -\frac{1}{2}\overline B^T\overline BD_ze + \frac{1}{2}\overline B^T v \label{eq: zl} \\
	\dot e &= -D_z^T\overline B^T\overline BD_ze + D_z^T\overline B^T v \label{eq: el} \\
	\dot v &= -\frac{1}{2}\overline BD_ze - v \label{eq: vl},
\end{align}
where $e$ is the stacked vector of $e_k$'s for all $k\in\{1,\dots,|\mathcal{E}|\}$. Define the speed of the agent $i$ by
\begin{equation*}
s_i \dfb ||v_i||,
\end{equation*}
whose time derivative is given by $	\dot s_i = \frac{v_i^T\dot{v_i}}{s_i}$. The compact form involving all the agents' speed can be written as
\begin{equation}
	\dot s = D_{\tilde s}D_v^T\dot v = -\frac{1}{2}D_{\tilde s}D_v^T\overline BD_ze - D_{\tilde s}D_v^Tv, \label{eq: sl}
\end{equation}
where $s$ and $\tilde s$ are the stacked vectors of $s_i$'s and $\frac{1}{s_i}$'s for all $i\in\{1,\dots,|\mathcal{V}|\}$ respectively. Now we are ready to show the local exponential convergence to the origin of the speed of the agents and the error distances in the edges.
\vspace{5pt}
\begin{lemma}
	\label{lem: expH}
	The origins $e = \mathbf{0}$ and $s = \mathbf{0}$ of the error and speed systems derived from (\ref{eq: H}) are locally exponentially stable if the given desired shape $\mathcal{D}$ is infinitesimally and minimally rigid.
\end{lemma}
\vspace{5pt}
\begin{IEEEproof}
Consider  $e_\lambda$ and $s_\lambda$ as the stacked vectors of the error signals $e_k$ and speeds $s_k$ derived from (\ref{eq: Hl}) for any $\lambda \in [0,1]$, which includes the system (\ref{eq: H}) for $\lambda = 1$. From the definition of $e_k$, we know that all the $e_\lambda$ share the same stability properties by invoking Lemma \ref{lem: H}, so do $s_\lambda$ as well.

Consider the following candidate Lyapunov function for the autonomous system (\ref{eq: zl})-(\ref{eq: sl}) derived from (\ref{eq: Hl}) with $\lambda = 0.5$
\begin{equation*}
	V = \frac{1}{2}||e||^2 + ||s||^2,
\end{equation*}
whose time derivative satisfies
\begin{align}
	\frac{\mathrm{d}V}{\mathrm{d}t} &= e^T\dot e + 2s^T\dot s \nonumber \\
&= -e^TD_z^T\overline B^T\overline BD_ze + e^TD_z^T\overline B^T v -
\underbrace{s^TD_{\tilde s}D_v^T}_{v^T}\overline BD_ze \nonumber \\ 
	&- 2\underbrace{s^TD_{\tilde s}}_{\mathbf{1}_{1\times|\mathcal{V}|}}D_v^T v \nonumber \\
	&= -e^TD_z^T\overline B^T\overline BD_ze - 2||s||^2.
\end{align}
We first note that the elements of the matrix $D_z^T\overline B^T\overline BD_z$ are of the form $z_i^Tz_j$ with $i, j\in\{1, \dots, |\mathcal{E}|\}$. It has been shown in \cite{MouMorseBelSunAnd15} that for minimally rigid shapes these dot products can be expressed as a (local) smooth functions of $e$, i.e. $z_i^Tz_j = g_{ij}(e)$, allowing us to write $D_z^T\overline B^T\overline BD_z = Q(e)$. For infinitesimally minimally rigid frameworks $R(z)$ is full rank, so we have that $Q(\mathbf{0}) = R(z^*)R^T(z^*)$ is positive definite with $z^*\in\mathcal{D}$. Moreover, since the eigenvalues of a matrix are continous functions of their entries, we have that $Q(e)$ is positive definite in the compact set $\mathcal{Q}\dfb \{e: ||e||^2 \leq \rho\}$ for some $\rho > 0$. Therefore, if the initial conditions for the error signal and the speed satisfy $||e(0)||^2 + ||s(0)||^2 \leq \rho$, since $V$ is not increasing we have that
\begin{equation}
	\frac{\mathrm{d}V}{\mathrm{d}t} \leq -\sigma_{min}||e||^2 - 2||s||^2,
\end{equation}
	where $\sigma_{\text{min}} > 0$ is the smallest (and always positive) eigenvalue of $Q(e)$ in the compact set $\mathcal{Q}$. Hence we arrive at the local exponential convergence of $e(t)$ and $s(t)$ to the origin.
\end{IEEEproof}
\begin{remark}
	It is worth noting that the region of attraction determined by $\rho$ in the proof of Lemma \ref{lem: expH} for $\lambda = 0.5$ might be different from the one for $\lambda = 1$, since Lemma \ref{lem: H} only refers to the Jacobian of (\ref{eq: Hl}), i.e., the linearization of the system about the equilibrium.
\end{remark}
It can be concluded from the exponential convergence to zero of the speeds of the agents $s(t)$ that the formation will eventually stop. This implies that $p(t)$ will converge exponentially to a finite point in $\R^m$ as $z(t)$ converges exponentially to $\mathcal{D}$.

\subsection{Robustness issues caused by mismatches}
\label{sec: is2}
Consider a distance-based formation control problem with $n=2$. It is not difficult to conclude that if the two agents do not share the same prescribed distance to maintain, then an eventual steady-state motion will happen regardless of the dynamics of the agents since the agent with a smaller prescribed distance will chase the other one. Therefore, for $n > 2$ it would not be surprising to observe some collective motion in the steady-state of the formation if the neighboring agents do not share the same prescribed distance to maintain.

Consider that two neighboring agents disagree on the desired squared distance $d_k^2$ in between, namely
\begin{equation}
d_k^{2\,\text{tail}} = d_k^{2\,\text{head}} - \mu_k,
\label{eq: mu}
\end{equation}
where $d_k^{\text{tail}}$ and $d_k^{\text{head}}$ are the different desired distances that the agents $i$ and $j$ respectively in $\mathcal{E}_k=(i,j)$ want to maintain for the same edge, so $\mu_k\in\R$ is a constant mismatch. It can be checked that this disagreement leads to \emph{mismatched potential functions}, therefore agents $i$ and $j$ do not share anymore the same $V_k$ for $\mathcal{E}_k = (i,j)$, namely
\begin{equation*}
	V^i_k = \frac{1}{4}(||z_k||^2 - d_k^2 + \mu_k)^2 \nonumber, \quad
	V^j_k = \frac{1}{4}(||z_k||^2 - d_k^2)^2 \nonumber,
\end{equation*}
under which the control laws for agents $i$ and $j$ use the gradients of $V^i_k$ and $V^j_k$ respectively for the edge $\mathcal{E}_k = (i,j)$. In the presence of one mismatch in every edge, the control signal (\ref{eq: un}) can be rewritten as
\begin{equation}
	u = -v -\overline BD_ze - \overline S_1D_z\mu, \label{eq: umu}
\end{equation}
where $S_1$ is constructed from the incidence matrix by setting its $-1$ elements to $0$, and $\mu\in\R^{|\mathcal{E}|}$ is the stacked column vector of $\mu_k$'s for all $k\in\{1, \dots, |\mathcal{E}|\}$. Note that (\ref{eq: umu}) can be also written as
\begin{equation}
	u = -v -\overline BD_ze - \overline A_1(\mu)z, \label{eq: umu2}
\end{equation}
where the elements of $A_1$ are
\begin{equation}
	a_{ik} \dfb \begin{cases}\mu_k & \text{if} \quad i = {\mathcal{E}_k^{\text{tail}}} \\
	0 & \text{otherwise.}
	\end{cases}
	\label{eq: A1}
\end{equation}

Inspired by \cite{MouMorseBelSunAnd15}, we will show how $\mu$ can be seen as a parametric disturbance in an autonomous system whose origin is exponentially stable. Consider the dynamics of the error signal $e$ and the speed of the agents $s$ derived from system (\ref{eq: pdyn}) with the control input (\ref{eq: un})
\begin{align}
	\dot e &= 2D_z^T\overline B^T v \label{eq: e} \\
	\dot s &= -s - D_{\tilde s}D_v^T\overline BD_ze,
\end{align}
and define
\begin{align}
	\alpha_{ki} &= z_k^Tv_i, \quad k\in\{1,\dots,|\mathcal{E}|\}, i\in\{1,\dots,|\mathcal{V}|\} \label{eq: alpha} \\
	\beta_{ij} &= v_i^Tv_j, \quad i,j\in\{1,\dots,|\mathcal{V}|\}, i\neq j. \label{eq: beta}
\end{align}
We stack all the $\alpha_{ki}$'s and $\beta_{ij}$'s in the column vectors $\alpha\in\R^{|\mathcal{E}||\mathcal{V}|}$ and $\beta\in\R^{\frac{|\mathcal{V}|(|\mathcal{V}|-1)}{2}}$ respectively and define $\gamma \dfb \begin{bmatrix}e^T & s^T & \alpha^T & \beta^T \end{bmatrix}^T$. We recall the result from Lemma \ref{lem: expH} that for any infinitesimally and minimally rigid framework there exists a neighborhood $\mathcal{U}_z$ about this framework such that for all $z_i,z_j\in\mathcal{U}_z$ with $i,j\in\{1,\dots,|\mathcal{E}|\}$ we can write $z_i^Tz_j$ as a smooth function $g_{ij}(e)$. Then using (\ref{eq: e})-(\ref{eq: beta}) we get
\begin{equation}
	\dot\gamma = f(\gamma), \label{eq: f}
\end{equation}
which is an autonomous system whose origin is locally exponentially stable using the results from Lemmas \ref{lem: H} and \ref{lem: expH}. Obviously, in such a case, the following Jacobian evaluated at $\gamma = \mathbf{0}$
\begin{equation*}
	J = \left.\frac{\partial f(\gamma)}{\partial \gamma}\right|_{\gamma = \mathbf{0}},
\end{equation*}
has all its eigenvalues in the left half complex plane. From the system (\ref{eq: pdyn}) with control law (\ref{eq: umu}) we can \emph{extend} (\ref{eq: f}) but with a parametric disturbance $\mu$ because of the third term in (\ref{eq: umu}), namely
\begin{equation}
	\dot\gamma = f(\gamma, \mu), \label{eq: fmu}
\end{equation}
where $f(\gamma, \mathbf{0})$ is the same as in (\ref{eq: f}) derived from the gradient controller. Therefore, for a sufficiently small $||\mu||$, the Jacobian $\left.\frac{\partial f(\gamma, \mu)}{\partial \gamma}\right|_{\gamma = 0}$ is still a stable matrix since the eigenvalues of a matrix are continuous functions of its entries. Although system (\ref{eq: fmu}) is still stable under the presence of a small disturbance $\mu$, the equilibrium point is not the origin in general anymore but $\gamma(t) \to \hat \gamma(\mu)$ as $t$ goes to infinity, where $\gamma_\mu \dfb \hat \gamma(\mu)$ is a smooth function of $\mu$ with zero value if $\mu = \mathbf{0}$ \cite{khalil1996nonlinear}. This implies that in general each component of $e, s, \alpha$ and $\beta$ converges to a non-zero constant with the following two immediate consequences: the formation shape will be distorted, i.e., $e \neq \mathbf{0}$; and the agents will not remain stationary, i.e., $s \neq \mathbf{0}$. The meaning of having in general non-zero components in $\alpha$ is that the vector velocities of the agents have a fixed relation with the steady-state shape, while the fixed components in $\beta$ denote a constant relation between the vectors of the agents' velocities. If the disturbance $||\mu||$ is sufficiently small, then $||\gamma_\mu|| < \rho$ for some small $\rho\in\R^+$ implying that $||e_\mu|| < \rho$, and if further $\rho$ is sufficiently small, then the stationary distorted shape is also infinitesimally and minimally rigid. In addition, since the speeds of the agents converge to a constant (in general non-zero constant), then only translations and/or rotations of the stationary distorted shape can happen.
\begin{theorem}
	Consider system (\ref{eq: pdyn}) with the \emph{mismatched} control input (\ref{eq: umu}) where the desired shape for the formation is infinitesimally and minimally rigid. There exist sufficiently small $\epsilon_1,\epsilon_2\in\R^+$ such that if the mismatches satisfy $||\mu|| \leq \epsilon_1$, then the error signal $e(t) \to e^*\in\R^{|\mathcal{E}|}$ as $t\to\infty$ with $||e^*|| \leq \epsilon_2$ such that the steady-state shape of the formation is still infinitesimally and minimally rigid but distorted with respect to the desired one. Moreover, the velocity of the agents $v_i(t)\to v_i^*(t)$ as $t\to\infty$, where all the $v_i^*(t), \forall i\in\{1,\dots,n\}$ define a steady-state collective motion that can be captured by constants angular and translational velocities $^b\omega^*$ and $^bv_c^*$, respectively, where $O_b$ has been placed in the centroid of the resultant distorted infinitesimally and minimally rigid shape.
\end{theorem}
\begin{IEEEproof}
	We have that system (\ref{eq: fmu}), derived from (\ref{eq: pdyn}) and (\ref{eq: umu}), is self-contained and its origin is locally exponentially stable with $\mu = \mathbf{0}_{|\mathcal{E}|\times 1}$. Then, a small parametric perturbation $\mu$ such that $||\mu||\leq \epsilon_1$, for some positive small $\epsilon_1$, does not change the exponential stability property of (\ref{eq: fmu}). However, the equilibrium point of (\ref{eq: fmu}) at the origin can be shifted. In particular, the new shifted equilibrium is a continuous function of $\mu$, therefore $e(t)\to e_\mu\in\R^{|\mathcal{E}|}$ as $t$ goes to infinity, where if $\epsilon_1$ is sufficiently small, then $||e_\mu|| \leq \epsilon_2$ such that the stationary shape is infinitesimally rigid. We also have that the elements of $s(t)\to s_\mu$ as $t$ goes to infinity with $s_\mu\in\R^{|\mathcal{V}|}$. Note that the elements of $s_\mu$ are non-negative and in general non-zero. Hence, the agents will not stop moving in the steady-state. Since the steady-state shape of the formation locally converges to an infinitesimally and minimally rigid one, from the error dynamics (\ref{eq: e}) we have that
	\begin{equation*}
		D_{z(t)}^T\overline B^T v(t) = R\left(z(t)\right)v(t) \to \mathbf{0}_{m|\mathcal{V}|\times 1}, \quad t\to\infty,
	\end{equation*}
	therefore $v(t)\to v_\mu(t)$ as $t$ goes to infinity, where the non-constant $v_\mu(t)\in\R^{|\mathcal{V}|}$ belongs to the null space of $R\left(z_\mu(t)\right)$, $z_\mu(t)\in\mathcal{Z}_\mu$ and the set $\mathcal{Z}_\mu$ is defined as in (\ref{eq: Zp}) but corresponding to the inter-agent distances of the distorted infinitesimally and minimally rigid shape with $e=e_\mu$. Note that obviously, the evolution of $z(t)$ is a consequence of the evolution of agents' velocities in $v(t)$. The null space of $R\left(z_\mu\right)$ corresponds to the infinitesimal motions $\delta p_i$ for all $i$ such that all the inter-agent distances of the distorted formation are constant, namely
	\begin{equation*}
		R(z_\mu)\delta p = R(z_\mu)v_\mu\,\delta t = \mathbf{0}_{m|\mathcal{V}|\times 1},
	\end{equation*}
or in order words
\begin{equation}
	v_i(t)\to v_{\mu_i}(t), \quad t\to\infty,
\end{equation}
where the velocities $v_{\mu_i}(t)$'s for all the agents are the result of rotating and translating the steady-state distorted shape defined by $\mathcal{Z}_\mu$. This steady-state collective motion of the distorted formation can be represented by the rotational and translational velocities $^b\omega^*(t)\in\R^m$ and $^bv_c(t)^*\in\R^m$ at the centroid of the distorted rigid shape. 

Now we are going to show that the velocities $^b\omega^*(t)$ and $^bv_c(t)^*$ are indeed constant. By definition we have that $||v_{\mu_i}(t)|| = s_{\mu_i}$. Since the speed $s_{\mu_i}$ for agent $i$ is constant but not its velocity $v_{\mu_i}(t)$, the acceleration $a_{\mu_i}(t) = \frac{\mathrm{d}v_{\mu_i}(t)}{\mathrm{d}t}$ is perpendicular to $v_{\mu_i}(t)$. The expression of $a_{\mu_i}(t)$ can be derived from (\ref{eq: umu2}) and is given by
\begin{equation}
	a_{\mu_i}(t) = -v_{\mu_i}(t) - \sum_{k=1}^{|\mathcal{E}|}b_{ik}z_{\mu_k}(t)e_{\mu_k} +  \sum_{k=1}^{|\mathcal{E}|}a_{ik}z_{\mu_k}(t), \label{eq: amu}
\end{equation}
where $b_{ik}$ are the elements of the incidence matrix, and $a_{ik}$ are the elements of the \emph{perturbation} matrix $A_1$ as defined in (\ref{eq: A1}). From (\ref{eq: amu}) it is clear that the norm $||a_{\mu_i}(t)|| = \Gamma_i(\gamma_\mu)$ is constant. In addition, since $a_{\mu_i}(t)$ is a continuous function, i.e., the acceleration vector cannot switch its direction, and it is perpendicular to $v_{\mu_i}(t)$, the only possibility for the distorted formation is to follow a motion described by constant velocities $^b\omega^*$ and $^bv_c^*$ at its centroid.
\end{IEEEproof}
\begin{remark}
	In particular, in 2D the distorted formation will follow a closed orbit if $\Gamma_i(\gamma_\mu) \neq 0$ for all $i$, or a constant drift if $\Gamma_i(\gamma_\mu) = 0$ for all $i$. This is due to the fact that in 2D, $^b\omega^*$ and $^bv^*$ are always perpendicular or equivalently $a_{\mu_i}(t)$ and $v_{\mu_i}(t)$ lie in the same plane. The resultant motion in 3D is the composition of a drift plus a closed orbit, since $^b\omega^*$ and $^bv_c^*$ are constant and they do not need to be perpendicular to each other as it can be noted in Figure \ref{fig: motion}.
\end{remark}
\begin{remark}
	Although the disturbance $\mu$ acts on the acceleration of the second-order  agents, it turns out that the resultant collective motion has the same behavior as for having the disturbance $\mu$ acting on the velocity for first-order agents. A detailed description of such a motion related to the disturbance in first-order agents can be found in \cite{SunMouAndMor14,MaJaCa15}.
\end{remark}
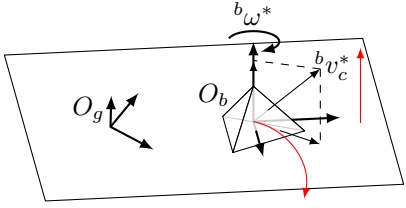
\begin{figure}
\centering
\input{FIGmotionTetra.tex}
\caption{The velocities $^b\omega^*$ and $^bv_c^*$ at the centroid of the tetrahedron rotates and translates the infinitesimally and minimally rigid formation respectively. If the velocitiy vectors $^b\omega^*$ and $^bv_c^*$ are constant, then the formation describes a closed orbit in the plane where the projection of $^bv_c^*$ over such a plane and $^b\omega^*$ is perpendicular (always the case in 2D formations) plus a constant drift along the projection of $^bv_c^*$ over $^b\omega^*$. Note that for the case where $^bv_c^*$ and $^b\omega^*$ are only parallel or perpendicular, then only the drift or the closed orbit motion, respectively, will occur.}
\label{fig: motion}
\end{figure}

\section{Estimator-based gradient control}
\label{sec: esti}
It is obvious that if for the edge $\mathcal{E}_k=(i,j)$ only one of the agents controls the desired inter-agent distance, then a mismatch $\mu_k$ cannot be present. However, this solution leads to a directed graph in the sensing topology and the stability of the formation can be compromised \cite{BaArWe11}. It is desirable to maintain the undirected nature of the sensing topology since it comes with intrinsic stability properties. Then, the control law (\ref{eq: umu}) must be \emph{augmented} in order to remove the undesired effects described in Section \ref{sec: issues}. A solution was proposed in \cite{MarCaoJa15} for first-order agents consisting of estimators based on the internal model principle. For each edge $\mathcal{E}_k=(i,j)$, there is only one agent that is assigned to be the \emph{estimating agent} which is responsible for running an estimator to calculate and compensate the associated $\mu_k$. The estimator proposed in \cite{MarCaoJa15} is conservative since the estimator gain has to satisfy a lower-bound (which can be explicitly computed based on the initial conditions) in order to guarantee the exponential stability of the system. Using such distributed estimators, all the undesirable effects are removed at the same time as the estimating agent calibrates its measurements with respect to the non-estimating agent. Another minor issue in the solution of \cite{MarCaoJa15} is that the estimating agents cannot be chosen arbitrarily. Here we are going to present an estimator for second-order agents where the estimating agents and the estimator gain can be chosen arbitrarily (thus, removing the restrictive conditions in \cite{MarCaoJa15}). The solution removes the effect of the undesired collective motion but at the cost of not achieving accurately the desired shape $\mathcal{Z}$, where a bound on the norm of the signal error $e(t)$ for all time $t$, however, can be provided. Furthermore, we will show that for the particular cases of the triangle and tetrahedron, the proposed estimator achieves precisely the desired shapes.

Let us consider the following distributed control law with estimator
\begin{equation}
	\begin{cases}\dot{\hat\mu} &= \hat u \\
u &= -v -\overline BD_ze - \overline S_1D_z(\mu - \hat\mu),\end{cases} \nonumber
\end{equation}
where $\hat\mu\in\R^{|\mathcal{E}|}$ is the estimator state and $\hat u$ is the estimator input to be designed. Substituting the above control law to (\ref{eq: pdyn}) gives us the following autonomous system
\begin{align}
	\dot p &= v \label{eq: pest} \\
	\dot v &= -v -\overline BD_ze - \overline S_1D_z(\mu - \hat\mu) \label{eq: vest} \\
	\dot z &= \overline B^T \dot p = \overline B^T v \label{eq: zest} \\
	\dot e &= 2D_z^T \dot z = 2D_z^T\overline B^T v \label{eq: eest} \\
\dot{\hat\mu} &= \hat u. \label{eq: muhatest}
\end{align}
Note that the estimating agents are encoded in $S_1$, in other words, for the edge $\mathcal{E}_k$ the estimating agent is $\mathcal{E}_k^{\text{tail}}$. 
\begin{theorem}
	\label{th: est}
	Consider the autonomous system (\ref{eq: pest})-(\ref{eq: muhatest}) with non-zero mismatches and a desired infinitesimally and minimally rigid formation shape $\mathcal{Z}$. Consider also the following distributed control action for the estimator $\hat\mu$
\begin{equation}
\hat u = -D_z^T\overline S_1^Tv,\label{eq: uhat}
\end{equation}
	where the estimating agents are chosen arbitrarily. Then the equilibrium points $(p^*, v^*, z^*, e^*, \hat\mu^*)$ of (\ref{eq: pest})-(\ref{eq: muhatest}) are asymptotically stable where $v^*=\mathbf{0}$ and the steady-state deformation of the shape satisfies $||e^*||^2\leq 2||\mu - \hat\mu(0)||^2 + 2||v(0)||^2 + ||e(0)||^2$. Moreover, for the particular cases of triangles and tetrahedrons, the equilibrium $e^*=\mathbf{0}$, i.e. $\hat\mu^* = \mu$ and $z^*\in\mathcal{Z}$.
\end{theorem}
\begin{IEEEproof}
	First we start proving that (\ref{eq: uhat}) is a distributed control law. This is clear since the dynamics of $\hat\mu_k$ (the $k$'th element of $\hat\mu$) are given by
	\begin{equation}
		\dot{\hat\mu}_k = z_k^Tv_{\mathcal{E}_k^{\text{tail}}}, \label{eq: Vest}
	\end{equation}
which implies that the estimating agent $\mathcal{E}_k^{\text{tail}}$ for the edge $\mathcal{E}_k$ is only using the dot product of the associated relative position $z_k$ and its own velocity. Note that using the notation in (\ref{eq: alpha}), the above estimator input is given by $\alpha_{k\mathcal{E}_k^{\text{tail}}}$. Consider the following Lyapunov function candidate
\begin{equation*}
V = \frac{1}{2}||\xi||^2 + \frac{1}{2}||v||^2 + \frac{1}{4}||e||^2,
\end{equation*}
with $\xi = \mu - \hat\mu$, which satisfies
\begin{align}
	\frac{\mathrm{d}V}{\mathrm{d}t} &= \xi^T\dot\xi + v^T\dot v + \frac{1}{2}e^T\dot e \nonumber \\
									&= \xi^TD_z^T\overline S_1^Tv - ||v||^2 - v^T\overline BD_ze - v^T\overline S_1D_z\xi \nonumber \\ & + e^TD_z^T\overline B^Tv \nonumber \\
				  &= -||v||^2. \label{eq: Vdest}
\end{align}
From this equality we can conclude that $\xi, v$ and $e$ are bounded. Moreover, from the definition of $e$, $z$ is also bounded. Thus all the states of the autonomous system (\ref{eq: vest})-(\ref{eq: muhatest}) are bounded, so one can conclude the convergence of $v(t)$ to zero in view of (\ref{eq: Vdest}). Furthermore, since the right-hand side of (\ref{eq: vest}) is uniformly continuous, $\dot v(t)$ converges also to zero by Barbalat's lemma. By invoking the LaSalle's invariance principle, looking at (\ref{eq: vest}) the states $e, \xi$ and $z$ converge asymptotically to the largest invariance set given by
\begin{equation}
	\mathcal{T} \dfb \{e, z, \xi : \overline S_1D_z\xi + \overline BD_ze = \mathbf{0}_{m|\mathcal{V}|\times 1}\}, \label{eq: Tp}
\end{equation}
in the compact set
\begin{equation}
	\mathcal{Q}\dfb \{\xi,v,e : ||\xi||^2 + ||v||^2 + \frac{1}{2}||e||^2 \leq \rho\}, \label{eq: Qr}
\end{equation}
with $0<\rho\leq 2V(0)$. Since $v=\mathbf{0}$ for all points in this invariant set, it follows from (\ref{eq: zest})-(\ref{eq: muhatest}) that $z, e$ and $\hat\mu$ are constant in this invariant set. In other words, $z(t)\to z^*, e(t)\to e^*$ and $\xi(t)\to \xi^*$ as $t$ goes to infinity, where $z^*, e^*$ and $\hat\mu^*$ are fixed points satisfying (\ref{eq: Tp}). Note that by comparing (\ref{eq: pest}) and (\ref{eq: zest}) one can also conclude that $p(t) \to p^*$ as $t$ goes to infinity. In general we have that $e^*$ and $\xi^*$ are not  zero vectors, therefore $z^* \notin\mathcal{Z}$. It is also clear that $||e^*||^2\leq 2\rho$, therefore for a sufficiently small $\rho$, the resultant (distorted) formation will also be infinitesimally and minimally rigid.

Now we are going to show that $e^*, \xi^* = \mathbf{0}$ for triangles and tetrahedrons. Since triangles and tetrahedrons are derived from complete graphs, the distorted shape when $\rho$ is sufficiently small will also be a triangle or a tetrahedron, i.e., we are excluding non-generic situations, e.g., collinear or coplanar alignments of the agents in $\R^2$ or $\R^3$. In the triangular case we have two possibilities after choosing the estimating agents: their associated directed graph is cyclic (each agent estimates one mismatch) or acyclic (one agent estimates two mismatches and one of the other two agents estimate the remaining mismatch).

The cyclic case for the estimating agents in the triangle corresponds to the following matrices
\begin{equation*}
B = \begin{bmatrix}1 & 0 & -1 \\ -1 & 1 & 0 \\ 0 & -1 & 1\end{bmatrix}, \,
S_1 = \begin{bmatrix}1 & 0 & 0 \\ 0 & 1 & 0 \\ 0 & 0 & 1 \end{bmatrix},
\end{equation*}
and by substituting them into the equilibrium condition in $\mathcal{T}$ we have that
\begin{equation}
	\left.
	\begin{array}{rc}
		z_1^*e^*_1 - z_3^*e^*_3 + z^*_1\xi_1^* &= 0 \\
		z_2^*e^*_2 - z_1^*e^*_1 + z^*_2\xi_2^* &= 0 \\
z_3^*e^*_3 - z_2^*e^*_2 + z^*_3\xi_3^* &= 0 \end{array}\right\}.
		\label{eq: S1I}
\end{equation}
Since the stationary formation is also a triangle for a sufficiently small $\rho$, then $z_1^*, z_2^*$ and $z_3^*$ are linearly independent. Therefore from (\ref{eq: S1I}) we have that $e_3^*, e_1^*, e_2^* = 0$ respectively and consequently we have that $\xi_1^*, \xi_2^*, \xi_3^* = 0$. 

Without loss of generality the acyclic case for the estimating agents in the triangle corresponds to the following matrices
\begin{equation}
B = \begin{bmatrix}-1 & 0 & -1 \\ 1 & 1 & 0 \\ 0 & -1 & 1\end{bmatrix}, \,
S_1 = \begin{bmatrix}0 & 0 & 0 \\ 1 & 1 & 0 \\ 0 & 0 & 1 \end{bmatrix},
	\label{eq: Bt}
\end{equation}
and by substituting them into the equilibrium condition in $\mathcal{T}$ we have that
\begin{equation}
	\left.
	\begin{array}{rc}
-z_1^*e^*_1 - z_3^*e^*_3 &= 0 \\
z_2^*e^*_2 + z_1^*e^*_1 + z^*_2\xi_2^* + z^*_1\xi_1^* &= 0 \\
z_3^*e^*_3 - z_2^*e^*_2 + z^*_3\xi_3^* &= 0 \end{array}\right\}.
		\label{eq: S12}
\end{equation}
It is immediate from the first equation in (\ref{eq: S12}) that $e_1^*, e_3^* = 0$ and then from the third equation in (\ref{eq: S12}) we derive that $e_2^*, \xi_3^* = 0$, and hence $\xi_1^*, \xi_2^* = 0$ from the second equation in (\ref{eq: S12}).

For the sake of brevity we omit the proof for the tetrahedrons, but analogous to the analysis for the triangles, the key idea behind the proof is that the three relative vectors associated to an agent are linearly independent in 3D.
\end{IEEEproof}
\begin{remark}
	The solution proposed in Theorem \ref{th: est} is distributed in the sense that each agent runs its own estimator based on only local information in order to address or solve the \emph{global} problem of a distorted and moving formation.
\end{remark}
\begin{remark}
	Since the estimating agents are defined by $S_1$, by \emph{arbitrarily chosen} we mean that for a given edge in the incidence matrix $B$, the order of the agents (tail and head) can be chosen arbitrarily.
\end{remark}

\begin{remark}
	A wide range of shapes can be constructed through piecing together triangles or tetrahedrons. For example, a multi-agent deployment of an arbitrary geometric shape in 2D can be realized by employing a mesh of triangles. If we consider the case when the interaction between different sub-groups of triangular formations defines a directed graph (so without mismatches), but each groups' graphs are undirected triangles or tetrahedrons, then one can employ the results of Theorem \ref{th: est}.
\end{remark}

The use of distributed estimators in Theorem \ref{th: est}, except for triangles and tetrahedrons, does not prevent the undesirable effect of having a distortion in the steady-state shape with respect to the desired one. Nevertheless, the error norm $||e||$ is bounded by a constant $\sqrt{2\rho} > 0$. Since $||e||$ is a combination of all errors in every edge, we cannot use $\rho$ to prescribe a zero asymptotic error for some focus edges or to concentrate the bound only on some edges. This property is relevant if we want to reach a prescribed distance for a high-degree of accuracy for some edges. By exploiting the result in Theorem \ref{th: est} for triangles and tetrahedrons, we can construct a network topology (based on a star topology) that enables us to impose the error bound only on one edge while guaranteeing that the other errors converge to zero. We show this in the following proposition.
\begin{proposition}
	\label{pr: star}
	Consider the same mismatched formation control system as in Theorem \ref{th: est} with the equilibrium set given by $v = \mathbf{0}_{m|\mathcal{V}|\times 1}$ and $\mathcal{T}$ as in (\ref{eq: Tp}). Consider the triangular formation defined by $B$ and $S_1$ as in (\ref{eq: Bt}) for the incidence matrix and the estimating agents respectively. For any new agent $i, i \geq 4$ added to the formation, if we only link it to the agents $2$ and $3$, and at the same time we let the agents $2$ and $3$ to be the estimating agents for the mismatches in the new added links, then for a sufficiently small $\rho$ as in (\ref{eq: Qr})
\begin{equation}
	\lim_{t\to\infty}e_k(t) = 0, \, \forall k \neq 2,
\end{equation}
and $|e_2(t)|\leq\sqrt{2\rho}$ for all $t$.
\end{proposition}
\begin{IEEEproof}
	Clearly a star topology has been used for the newly added agents, where the center is the triangle formed by agents $1, 2$ and $3$. Note that the newly added agent $i, i\geq 4$ is forming a triangle with agents $2$ and $3$. Therefore, as explained in Theorem \ref{th: est}, if $\rho$ is sufficiently small, then the resultant distorted formation is also formed by triangles. We prove the claim by induction. First we derive the equations from $\mathcal{T}$ as in (\ref{eq: T}) for the proposed star topology with four agents
\begin{equation}
\left.
\begin{array}{rc}
-z_1^*e^*_1 - z_3^*e^*_3 &= 0 \\
z_2^*e^*_2 + z_1^*e^*_1 - z^*_4e_4 + z^*_2\xi_2^* + z^*_1\xi_1^* + z^*_4\xi_4&= 0 \\
z_3^*e^*_3 - z_2^*e^*_2 - z^*_5e_5 + z^*_3\xi_3^* + z^*_5\xi_5 &= 0  \\ 
z_4^*e^*_4 + z^*_5e^*_5 &= 0\end{array}\right\}.
	\label{eq: S123}
\end{equation}
As explained in the last part of the proof of Theorem \ref{th: est}, it is clear that the errors $e_1^*, e_3^*, e_4^*$ and $e_5^*$ must be zero and $e_2^2\leq\rho$. For any newly added agent $i\geq5$, we add a new equation to (\ref{eq: S123}) of the form
\begin{equation}
	z_l^*e^*_l + z_{l+1}^*e^*_{l+1} = 0, \label{eq: ll}
\end{equation}
where $l$ and $l+1$ are the labels of the two newly added edges. Thus for a sufficiently small $\rho$ we have that $z^*_l$ and $z^*_{l+1}$ are linearly independent so $e^*_l, e^*_{l+1} = 0$.
\end{IEEEproof}

It is possible to be more accurate in the estimation of the mismatches under mild conditions. We in \cite{MarCaoJa15} have proposed the following control law for the estimators in order to remove effectively both, the distortion and the steady-state collective motion
\begin{equation}
	\hat u_k = \kappa(e_k + \mu_k - \hat\mu_k), k\in\{1, \dots, |\mathcal{E}|\},
	\label{eq: uk}
\end{equation}
where $\kappa \in\R^+$ is a sufficiently high gain to be determined. Consider the following change of coordinates $h_k = e_k + \mu_k - \hat\mu_k$ and let $h\in\R^{|\mathcal{E}|}$ be the stacked vector of $h_k$'s for all $k\in\{1,\dots,|\mathcal{E}|\}$. By defining $S_2 \dfb B - S_1$ it can be checked that the following autonomous system derived from (\ref{eq: vest})-(\ref{eq: muhatest})
\begin{align}
	\dot v &= -v -\overline S_2D_ze - \overline S_1D_z h \label{eq: va}\\
	\dot e &= 2D_z\overline B^Tv \label{eq: ea} \\
	\dot h  &= 2D_z\overline B^Tv - \kappa h \label{eq: aa} \\
	\dot z &= \overline B^T v, \label{eq: za}
\end{align}
has an equilibrium at $e= h = \mathbf{0}$, $v =\mathbf{0}$ and $z^*\in\mathcal{Z}$. The linearization of the autonomous system (\ref{eq: va})-(\ref{eq: za}) about such an equilibrium point leads to
\begin{align}
\begin{bmatrix}\dot v \\ \dot e \\ \dot h \\ \dot z\end{bmatrix}
	=\begin{bmatrix}-\overline I_{|\mathcal{V}|} & -\overline S_2D_{z^*} & -\overline S_1D_{z^*} & \mathbf{0}  \\ 2D_{z^*}\overline B^T & \mathbf{0} & \mathbf{0} & \mathbf{0} \\
	2D_{z^*}\overline B^T & \mathbf{0} & -\kappa I_{|\mathcal{E}|} & \mathbf{0} \\
	\overline B^T & \mathbf{0} & \mathbf{0} & \mathbf{0}
\end{bmatrix}\begin{bmatrix}v \\ e \\ h \\ z\end{bmatrix}.
	\label{eq: lin}
\end{align}
From the Jacobian in (\ref{eq: lin}) we know that the stability of the system only depends on $v, e$ and $\alpha$ and not on $z$. We consider the following assumption as in \cite{MarCaoJa15}.
\vspace{5px}
\begin{assumption}
	The matrix $F\dfb \left[\begin{smallmatrix} -\overline I_{|\mathcal{V}|}& -\overline S_2D_{z^*} \\
2D_{z^*}\overline B^T & \mathbf{0}\end{smallmatrix}\right]$ is Hurwitz.
	\label{as: 1}
\end{assumption}
\begin{theorem}
	\label{th: est2}
	There exists a positive constant $\kappa^*$ such that the equilibrium corresponding to $\tilde\mu=\mu, v=\mathbf{0}$ and $e=\mathbf{0}$ (with $z^*\in\mathcal{Z}$) of the autonomous system (\ref{eq: vest})-(\ref{eq: muhatest}) with the estimator law (\ref{eq: uk}) is locally exponentially stable under Assumption \ref{as: 1} for all $\kappa > \kappa^* > 0$.
\end{theorem}
\begin{IEEEproof}
The resulting Jacobian by the linearization of (\ref{eq: va})-(\ref{eq: za}), derived from (\ref{eq: vest})-(\ref{eq: muhatest}), evaluated at the desired shape $z^*\in\mathcal{Z}$ is given by (\ref{eq: lin}). Note that the stability of the linear system (\ref{eq: lin}) does not depend on $z$, therefore the (marginal) stability of (\ref{eq: lin}) is given by analyzing the eigenvalues of the first $3\times 3$ blocks, i.e., the dynamics of $v, e$ and $h$. Also note that by Assumption \ref{as: 1} the first $2\times 2$ blocks, i.e., the matrix $F$, of the matrix in (\ref{eq: lin}) is Hurwitz. In addition, one can also easily check that the third block, i.e., $-\kappa I_{|\mathcal{E}|}$, of the main diagonal of (\ref{eq: lin}) is negative definite. We consider the following Lyapunov function candidate
	\begin{equation}
		V(v, e, h) = \begin{bmatrix}v^T & e^T\end{bmatrix}P\begin{bmatrix}v \\ e
	\end{bmatrix} + \frac{1}{2}||h||^2, \label{eq: Vh2}
	\end{equation}
where $P$ is a positive definite matrix such that $PF^T + FP = -2I$. For brevity, following the same arguments used in the proof of the main theorem in \cite{MarCaoJa15}, the time derivate of $V$ in a neighborhood of the equilibrium $h=\mathbf{0}, v=\mathbf{0}$ and $e=\mathbf{0}$ (with $z^*\in\mathcal{Z}$) can be given by
\begin{equation}
	\frac{\mathrm{d}V}{\mathrm{dt}} \leq -||v||^2 -||e||^2 + (m - \kappa)||h||^2,
\end{equation}
	where $m\in\R$ comes from the cross-terms that eventually are going to be dominated by $\kappa > \kappa^* = m$. We also refer to the book \cite{isidori2013nonlinear} for more details about the employed technique of requiring $F$ to be Hurwitz and employing the Lyapunov function (\ref{eq: Vh2}). This completes the proof.
\end{IEEEproof}
\begin{remark}
	Since $v(t)$ converges exponentially to zero, it follows immediately that $p(t)$ converges exponentially to a fixed point $p^*$.
\end{remark}
Assumption \ref{as: 1} is also related to the stability of formation control systems whose graph $\mathbb{G}$ defines a directed sensing topology. In fact, it is straightforward to check that the matrix in Assumption \ref{as: 1} is the Jacobian matrix for $v$ and $e$ in a distance-based formation control system (without mismatches) with only directed edges in $\mathbb{G}$, i.e., at the equilibrium (desired shape) the linearized non-linear system is given by
\begin{equation}
	\begin{bmatrix}
	\dot v \\ \dot e
	\end{bmatrix} =
	F	\begin{bmatrix}
	v \\ e
	\end{bmatrix}.
\end{equation}
Therefore, in order to satisfy Assumption \ref{as: 1}, one has to choose the estimating agents with the same topology as in stable directed rigid formations, e.g., a spanning tree. This selection can be checked in more detail in \cite{MarCaoJa15,4522613}.

These two presented strategies for the estimators in Theorems \ref{th: est} and \ref{th: est2} have advantages and drawbacks: the main advantage of (\ref{eq: uhat}) over (\ref{eq: uk}) is that we do not need to compute any gain and that there is a free choice of the estimating agents; on the other hand, (\ref{eq: uk}) guarantees exponential convergence to $\mathcal{Z}$ and a distributed estimation of $\mu$. It is also worth noting that a minor modification in (\ref{eq: uk}) also compensates time varying mismatches as it was studied in \cite{MarCaoJa15}.

\section{Motion control of second-order rigid formations}
\label{sec: motion}
In this section we are going to extend the findings in \cite{MaJaCa15} on the formation-motion control from the single-integrator to the double integrator case. In particular, we consider how to design the desired constant velocities $^b\omega^*$ and $^bv_c^*$ as in Figure \ref{fig: motion} for an infinitesimally and minimally rigid formation, i.e., for $e = 0$. The main feature of this approach is that the desired motion of the shape is designed with respect to $O_b$ but not $O_g$. Note that the latter is the common approach in the literature \cite{BaArWe11,oh2015survey}.

Following a similar strategy as in \cite{MaJaCa15, MaCDC}, we solve the motion control of rigid formation problem for double-integrator agents by employing mismatches as design parameters. More precisely, we assign two motion parameters $\mu_k$ and $\tilde\mu_k$ to agents $i$ and $j$ in the edge $\mathcal{E}_k = (i,j)$ resulting in the following control law\footnote{As a comparison, in \cite{MaJaCa15} the proposed controller is of the form of $u = -c_2\overline BD_ze + \overline A(\mu, \tilde\mu)z$.}
\begin{equation}
	u = -c_1 v -c_2\overline BD_ze + \overline A(\mu, \tilde\mu)z,
	\label{eq: uA}
\end{equation}
where $c_1, c_2\in\R^+$ are gains, $\mu\in\R^{|\mathcal{E}|}$ and $\tilde\mu\in\R^{|\mathcal{E}|}$ are the stacked vectors for all $\mu_k$ and $\tilde\mu_k$ and $A$ is defined by
\begin{equation}
	a_{ik} \dfb \begin{cases}\mu_k \quad \text{if} \quad i = {\mathcal{E}_k^{\text{tail}}} \\
		\tilde\mu_k \quad \text{if} \quad i = {\mathcal{E}_k^{\text{head}}} \\
						  0 \quad \text{otherwise.}
					\end{cases}
	\label{eq: A}
\end{equation}
The design of $^b\omega^*$ and $^bv_c^*$ is done via choosing appropriately the motion parameters $\mu$ and $\tilde\mu$, in the sense that we allow a desired steady-state collective motion but remove any distortion of the final shape. The design of the motion parameters in $A$ must take into account not only the desired acceleration but also the damping component in (\ref{eq: uA}) (which is different from the single-integrator case considered in \cite{MaJaCa15}).

Let the velocity error
\begin{equation}
	e_v = v - \overline A_v(\mu, \tilde\mu)z, \label{eq: Av}
\end{equation}
where $A_v(\mu, \tilde\mu)$ is designed employing the motion parameters described in \cite{MaJaCa15} directly related to the desired steady-state collective velocity. For the sake of completeness, we briefly describe how to compute $\mu$ and $\tilde\mu$ in $A_v$ for the prescribed $^bv_c^*$ and $^b\omega^*$ as in Figure \ref{fig: motion}. Since
\begin{equation}
	A_v(\mu, \tilde\mu) ^bz^* = \begin{bmatrix}\bar S_1 D_{^bz^*} & \bar S_2D_{^bz^*}\end{bmatrix}\begin{bmatrix}\mu \\ \tilde\mu\end{bmatrix} = T(^bz^*)\begin{bmatrix}\mu \\ \tilde\mu\end{bmatrix}, \label{eq: T}
\end{equation}
defines the steady-state velocity in the body frame, one can derive the following two conditions
\begin{align}
	\overline B^T\overline T(^bz^*)\begin{bmatrix}\mu \\ \tilde\mu\end{bmatrix} &= 0 \label{eq: Z} \\
		D_{^bz}\overline B^T\overline T(^bz^*)\begin{bmatrix}\mu \\ \tilde\mu\end{bmatrix} &= 0, \label{eq: E}
\end{align}
where (\ref{eq: Z}) stands for translations and (\ref{eq: E}) for rotations and translations, i.e., we set $\frac{\mathrm{d}z_k^*(t)}{\mathrm{dt}} = 0$ and $\frac{\mathrm{d}||z_k^*(t)||}{\mathrm{dt}} = 0$ for all $k\in\{1,\dots,|\mathcal{E}|\}$ in (\ref{eq: Z}) and (\ref{eq: E}) respectively. Let us split $\mu = \mu_v + \mu_\omega$ and $\tilde\mu = \tilde\mu_v + \tilde\mu_\omega$. In order to compute the distributed motion parameters $\mu_v, \tilde\mu_v$ for the translational velocity $^bv_c^*$ we eliminate the components of $\mu$ and $\tilde\mu$ that are not responsible for any motion by projecting the kernel of $\overline B^TT(^bz^*)$ onto the orthogonal space of the kernel of $T(^bz^*)$, namely
\begin{equation}
	\begin{bmatrix}\mu_v \\ \tilde\mu_v \end{bmatrix} \in \mathcal{\hat U}\dfb P_{\operatorname{Ker}\{T(^bz^*)\}^\bot}\left\{\operatorname{Ker}\{\overline B^TT(^bz^*)\} \right\}, \label{eq: U}
\end{equation}
where the operator $P_\mathcal{X}$ stands for the projection over the space $\mathcal{X}$. In a similar way, we need to remove the space responsible for the translational motion in the null space of the matrix in (\ref{eq: E}). Therefore, the computation of the distributed motion parameters $\mu_\omega, \tilde\mu_\omega$ for the rotational motion $^b\omega^*$ of the desired shape is obtained from (\ref{eq: E}) and (\ref{eq: U}) as
\begin{equation}
	\begin{bmatrix}\mu_\omega \\ \tilde\mu_\omega \end{bmatrix} \in  \mathcal{\hat W}\dfb P_{\mathcal{\hat U}^\bot}\left\{\operatorname{Ker}\{D_{^bz^*}^T\overline B^TT(^bz^*)\}\right\}.
\label{eq: W}
\end{equation}
When the velocity error $e_v$ is zero and we are at the desired shape $z^*(t)\in\mathcal{Z}$ with the desired velocities in $v^*(t)$, then from (\ref{eq: Av}) we have that
\begin{align}
	v^*(t) &= \overline A_v z^*(t) \\
	\dot v^*(t) &= \overline A_v \dot z^*(t) = \overline A_v\overline B^Tv^*(t) = \overline A_v\overline B^T \overline A_v z^*(t) \nonumber \\
	 &= \overline A_a z^*(t).
\end{align}
Note that the desired parameters in $A_a(\mu, \tilde\mu)$ correspond to the desired acceleration of the agents at the desired shape $\mathcal{Z}$. With this knowledge at hand, we can design the needed motion parameters for $A$ in the control law in (\ref{eq: uA}) as
\begin{equation}
A(\mu, \tilde\mu) = c_1A_v(\mu, \tilde\mu) + A_a(\mu, \tilde\mu),
\label{eq: A}
\end{equation}
since for $e_v = \mathbf{0}_{m|\mathcal{V}|\times 1}$ and $e = \mathbf{0}_{|\mathcal{E}|\times 1}$ the control law (\ref{eq: uA}) becomes
\begin{equation}
	u = \overline A_a(\mu, \tilde\mu)z^*(t). \label{eq: uAa}
\end{equation}
Note that $A(\mu, \tilde\mu)$ can be computed directly from the motion parameters for the desired velocity as in the first-order case. Therefore, there is no need of designing desired accelerations, which can be a more tedious task. We show in the following theorem that the desired collective-motion for the desired formation is stable for at least sufficiently small speeds.
\begin{theorem}
\label{th: rec}
	There exist constants $\rho, \rho_\mu, \epsilon, c_1, c_2>0$ for system (\ref{eq: pdyn}) with control law (\ref{eq: uA}), $A(\mu, \tilde\mu)$ as in (\ref{eq: A}) and with a given desired infinitesimally and minimally rigid shape $\mathcal{Z}$, such that if $[\begin{smallmatrix}\mu \\ \tilde\mu\end{smallmatrix}] \in \mathcal{M}\dfb\{\mu,\tilde\mu : ||[\begin{smallmatrix}\mu \\ \tilde\mu\end{smallmatrix}]||\leq\rho_\mu\}$, then the origin of the error dynamical system $e_v=\mathbf{0}$ and $e=\mathbf{0}$ corresponding to $z^*(t)\in\mathcal{Z}$ is exponentially stable in the compact set $\mathcal{Q}\dfb\{e, e_v : \frac{\epsilon c_1 + c_2}{4}||e||^2 + \frac{1}{2}||e_v||^2 \leq \rho\}$. In particular, the steady-state shape is the same as the desired one and the steady-state collective motion of the formation corresponds to $v^*(t) = \overline A_vz^*(t)$.
\end{theorem}
\begin{IEEEproof}
	First we rewrite the control law (\ref{eq: uA}) employing (\ref{eq: Av}) and (\ref{eq: A}) as
\begin{equation}
	u = -c_1e_v -c_2\overline BDze + \overline A_az. \label{eq: accU}
\end{equation}
Consider the following candidate Lyapunov function
\begin{equation}
	V = \frac{\epsilon c_1 + c_2}{4}||e||^2 + \frac{1}{2}||e_v||^2 + \epsilon e_v^T\overline BD_z e, \label{eq: Vmot}
\end{equation}
where $V$ is positive definite in a neighborhood about $e=\mathbf{0}$ and $e_v=\mathbf{0}$ for some sufficiently small $\epsilon\in\R^+$ in the compact set $\mathcal{Q}$ with $e = \mathbf{0}$ corresponding to $z\in\mathcal{Z}$. Note that even without knowing $c_1$ and $c_2$ yet, one can safely compute $\epsilon$ by assuming $c_1=0$ and $c_2=c_2^*$ for an arbitrary $c_2^*\in\R^+$. Then, later we restrict the choosing of $c_1$ and $c_2$ to be bigger than $0$ and $c_2^*$ respectively, since by this choice one does not change the positive semi-definite nature of (\ref{eq: Vmot}) for the calculated $\epsilon$ for $c_1 = 0$ and $c_2 = c_2^*$.

The time derivative of (\ref{eq: Vmot}) is given by
\small
\begin{align}
	\frac{\mathrm{d}V}{\mathrm{d}t} &= \frac{1}{2}(\epsilon c_1 + c_2)e^T\dot e + e_v^T\dot e_v + \epsilon e_v^T\overline BD_z \dot e +  \epsilon e^TD_z^T\overline B^T\dot e_v \nonumber \\ &+ \epsilon e_v^T\overline BD_{(\overline B^Tv)} e 
\nonumber \\
&= (\epsilon c_1 + c_2)e^TD_z^T\overline B^T(e_v+\overline A_v z) - c_1||e_v||^2-c_2e_v^T\overline BD_ze \nonumber \\
& + e_v^T\overline A_v\overline B^T\overline A_vz - e_v^T\overline A_v\overline B^T\overline A_vz - e_v^T\overline A_v\overline B^Te_v \nonumber \\ 
& + 2\epsilon e_v^T\overline BD_zD_z^T\overline B^Tv - c_1\epsilon e^TD_z^T\overline B^Te_v - c_2\epsilon e^TD_z^T\overline B^T\overline BD_ze \nonumber \\
&+ \epsilon e^TD_z^T\overline B^T\overline A_v\overline B^T\overline A_vz - \epsilon e^TD_z^T\overline B^T\overline A_v\overline B^T\overline A_vz \nonumber \\
&-\epsilon e^TD_z^T\overline B^T\overline A_v\overline B^T e_v + \epsilon e_v^T\overline BD_{(\overline B^Tv)} e \nonumber \\
&= (\epsilon c_1 + c_2)e^T\underbrace{D_z^T\overline B^T\overline A_v z}_{f_1(e,\mu,\tilde\mu)} - c_1||e_v||^2 - e_v^T\underbrace{\overline A_v\overline B^T}_{f_2(\mu,\tilde\mu)}e_v \nonumber \\
&+2\epsilon e_v^T\underbrace{\overline BD_zD_z^T\overline B^T}_{f_3(z)}e_v + 2\epsilon e_v^T\underbrace{\overline BD_zD_z^T\overline B^T\overline A_vz}_{f_4(\mu, \tilde\mu, z, e)} \nonumber \\
& - c_2\epsilon e^T\underbrace{D_z^T\overline B^T\overline BD_z}_{f_5(e)}e - \epsilon e^T\underbrace{D_z^T\overline B^T\overline A_v\overline B^T}_{f_6(\mu, \tilde\mu, z)} e_v \nonumber \\
& + \epsilon e_v^T\underbrace{\overline BD_{(\overline B^Tv)}}_{f_7(v)} e. \label{eq: Vdmot}
\end{align}
\normalsize
Since all the $f_i, i\in\{1, \dots, 7\}$ are locally Lipschitz functions in the compact sets $\mathcal{Q}$ and $\mathcal{M}$ and by using Young's inequality to every cross-term in (\ref{eq: Vdmot}), we can bound $\dot V$ as follows
\small
\begin{align}
	\frac{\mathrm{d}V}{\mathrm{d}t} &\leq \bigg(c_2\Big(M_1(\mu,\tilde\mu)-\epsilon\lambda_5\Big) + \epsilon c_1M_1(\mu,\tilde\mu) + \frac{3}{2}\bigg)||e||^2 \nonumber \\
 &+\bigg(-c_1 + M_2(\mu,\tilde\mu) \nonumber \\
 &+\epsilon^2 \Big(\frac{2\lambda_3}{\epsilon} +
M_4(\mu, \tilde\mu, z) + M_6(\mu, \tilde\mu, z) +
M_7\Big)\bigg)||e_v||^2,\label{eq: Vdmot2}
\end{align}
\normalsize
where $M_1$ and $M_4$ are related to the Lipschitz constant of $f_1$ and $f_4$
 in the compact set $\mathcal{Q}$ given $\mu$ and $\tilde\mu$, $M_2$ is the induced 2-norm of $f_2$ given $\mu$ and $\tilde\mu$, $M_6$ is the squared induced 2-norm of $f_4$ in the compact set $\mathcal{Q}$ given $\mu$ and $\tilde\mu$, and finally $M_7$ is the maximum squared induced 2-norm for $f_7$ in $\mathcal{Q}$ as well. We also have that $\lambda_3$ is the maximum eigenvalue of $f_3$ and $\lambda_5$ is the minimum eigenvalue of $f_5$ in the compact set $\mathcal{Q}$. First we note that for a sufficiently small $\rho$, $\lambda_5 > 0$ by the same argument of having a desired infinitesimally and minimally rigid formation as in Lemma \ref{lem: expH}. The time derivative (\ref{eq: Vdmot2}) can be made negative as a result of the following steps:
\begin{itemize}
	\item Choose a sufficiently small $\rho_\mu$ in $\mathcal{M}$ such that $M_1(\mu,\tilde\mu)-\epsilon\lambda_5<0$, i.e., downscale if necessary $\mu$ and $\tilde\mu$ by the same factor.
	\item Compute $M_2$ for the given $\mu,\tilde\mu$.
	\item Compute $M_4$ and $M_6$ for the given $\mu,\tilde\mu$ in the compact set $\mathcal{Q}$.
	\item Compute $M_7$ considering all the $v$ such that $||e_v||^2 \leq \rho$.
	\item Choose $c_1$ such that the second bracket in (\ref{eq: Vdmot2}) is negative.
	\item Given $c_1$ choose $c_2>c_2^*$ (employed for the calculation of $\epsilon$) such that the first bracket in (\ref{eq: Vdmot2}) is negative.
\end{itemize}
This guarantees the local exponential convergence of $e(t)$ and $e_v(t)$ to their origins, hence $z(t)\to z^*(t)\in\mathcal{Z}$, $v(t)\to \overline A_vz^*(t)$ and the stacked acceleration of the agents $a(t)\to \overline A_az^*(t)$ as $t$ goes to infinity.
\end{IEEEproof}
\begin{remark}
	The limitation given by $\rho_\mu$ is exclusively related to the desired speed of the agents\cite{MaJaCa15} and it does not restrict in any other way the desired collective motion for the formation. Therefore, once the motion parameters are given, for asserting the exponential stability of the system, one only has to downscale them if necessary. This (conservative) downscale has an intuitive physical explanation related to the condition of requiring $\lambda_5 > 0$ in Theorem \ref{th: rec}, e.g., the agents cannot be in a collinear configuration at any moment. In order to avoid such configurations where $\lambda_5$ is zero we need $e\in\mathcal{Q}$, i.e., for arbitrary velocities and positions of the agents such that $e, e_v \in \mathcal{Q}$, the third acceleration term in (\ref{eq: accU}) given by the motion parameters $\mu$ and $\tilde\mu$ must be small enough in order to avoid the possibility of pushing the agents to become collinear. Note that when the agents start close to the desired shape and have low velocity, the gradient terms in (\ref{eq: accU}) will not push the agents to such a configuration but to the desired shape. Nevertheless, imposing $\lambda_5 > 0$ for all time is a very conservative condition and we have verified in simulation that indeed it is not a necessary condition.
\end{remark}
\begin{remark}
	For $\mu,\tilde\mu = \mathbf{0}$ we have that $M_1, M_2, M_3, M_4, M_6 = 0$. Then employing (\ref{eq: Vmot}) and applying differently Young's inequalities in $(\ref{eq: Vdmot})$ one can prove that for $c_1, c_2 = 1$ the dissipative Hamiltonian system (\ref{eq: H}) is exponentially stable for a sufficiently small $\epsilon$.
\end{remark}
\begin{remark}
	For desired constant drifts in triangular and tetrahedron formations, it can be checked from \cite{MaJaCa15} that $M_1 = 0$. Therefore there is no restriction in the speed for such particular cases. In particular, one can use the Lyapunov function (\ref{eq: Vmot}) with $\epsilon = 0$. It turns out that the formation with the proposed motion-shape controller is asymptotically stable for any $c_2 > 0$ for $c_1 > ||\overline A_v(\mu,\tilde\mu) \overline B^T||^2$.
\end{remark}
\begin{remark}
	For a desired rotation about the centroid of an equilateral triangle, it can be checked from \cite{MaJaCa15} that in addition to $M_1 = 0$ we have that $e_v^T f_2(\mu,\tilde\mu)e_v = 0$ since $f_2$ is skew symmetric. Therefore by using (\ref{eq: Vmot}) with $\epsilon = 0$ one can prove the asymptotic stability of the origin of $e_v$ and $e$ for any $c_1, c_2 > 0$.
\end{remark}
\begin{remark}
	We recall the mismatched case in Section \ref{sec: is2}. So far, we can conclude that the upper bound for the distortion and speed of the resulting formation depends on the steady-state of the formed geometrical shape. This can be seen by looking at the third term on the right hand side of (\ref{eq: uA}) by just taking $\tilde\mu = \mathbf{0}$. On the other hand, there is a direct relation between a formation containing mismatches with a formation containing motion parameters by checking the two cases in (\ref{eq: umu}) and (\ref{eq: accU}). In the latter case we are able to calculate the steady-state of both shape and velocity. In the former (mismatched) case, we refer to the results in \cite{SunMouAndMor14} in order to check the approximated resultant motion. As a first step, the distortion in the mismatched case can be approached by comparing how far the mismatches are from the ones calculated employing (\ref{eq: U}), (\ref{eq: W}) and (\ref{eq: A}), where no distortion occurs.
\end{remark}

The result in Theorem \ref{th: rec} allows one to design the desired velocity for the given formation with respect to $O_b$ as in Figure \ref{fig: motion}. This result extends to the applications proposed in \cite{MaJaCa15} by using distributed motion parameters, such as steering an infinitesimally and minimally rigid formation by controlling the heading of the formation with respect to $O_g$, and for the tracking and enclosing of a target.

\section{Simulations}
\label{sec: exp}
In this section we validate the results in Theorems \ref{th: est}, \ref{th: est2} and \ref{th: rec} with numerical simulations.

\begin{figure}
\centering
\begin{subfigure}{0.48\columnwidth}
\includegraphics[width=1\columnwidth]{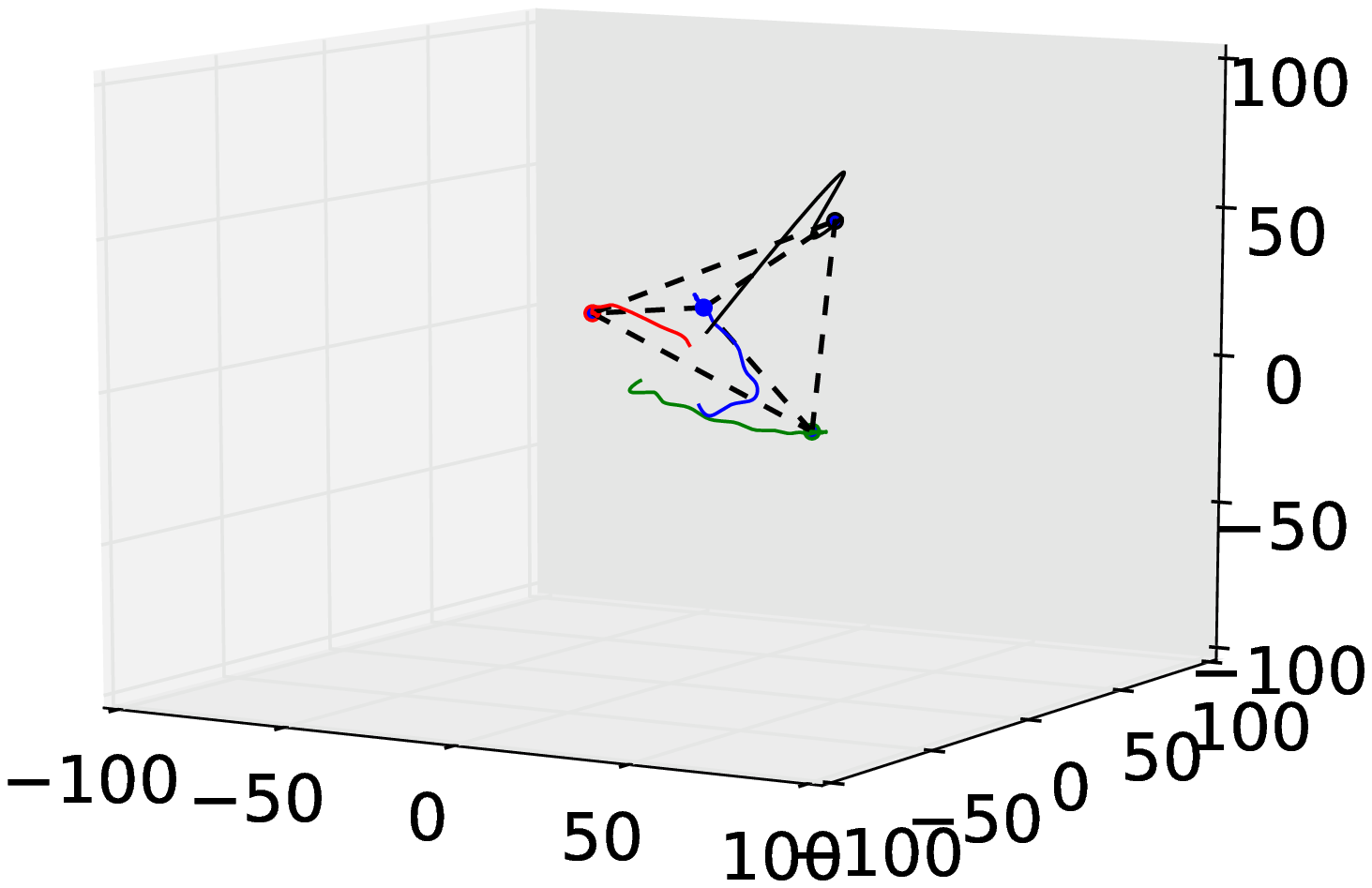}
\caption{}
\end{subfigure}
\begin{subfigure}{0.48\columnwidth}
\includegraphics[width=1\columnwidth]{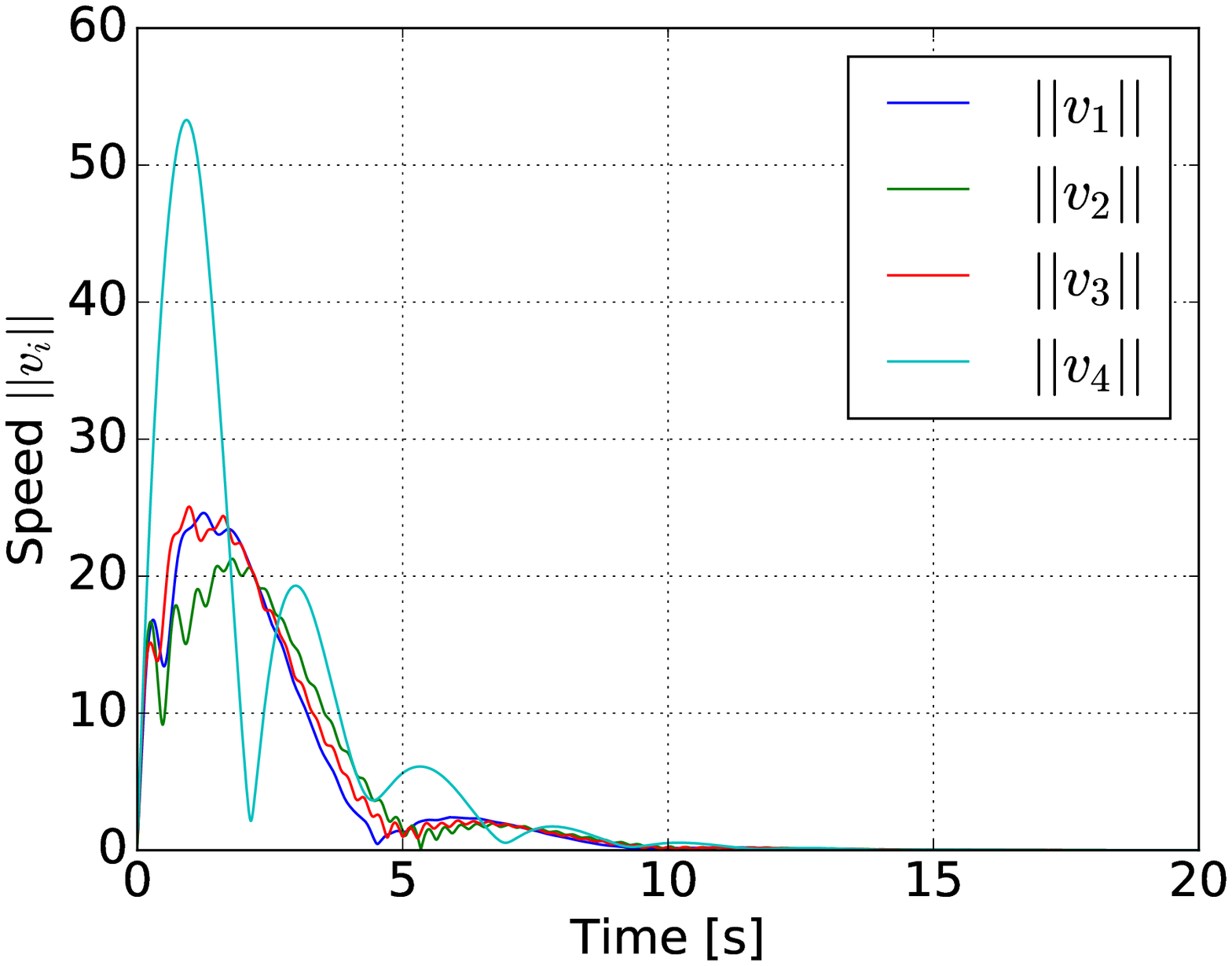}
\caption{}
\end{subfigure}
\begin{subfigure}{0.48\columnwidth}
\includegraphics[width=1\columnwidth]{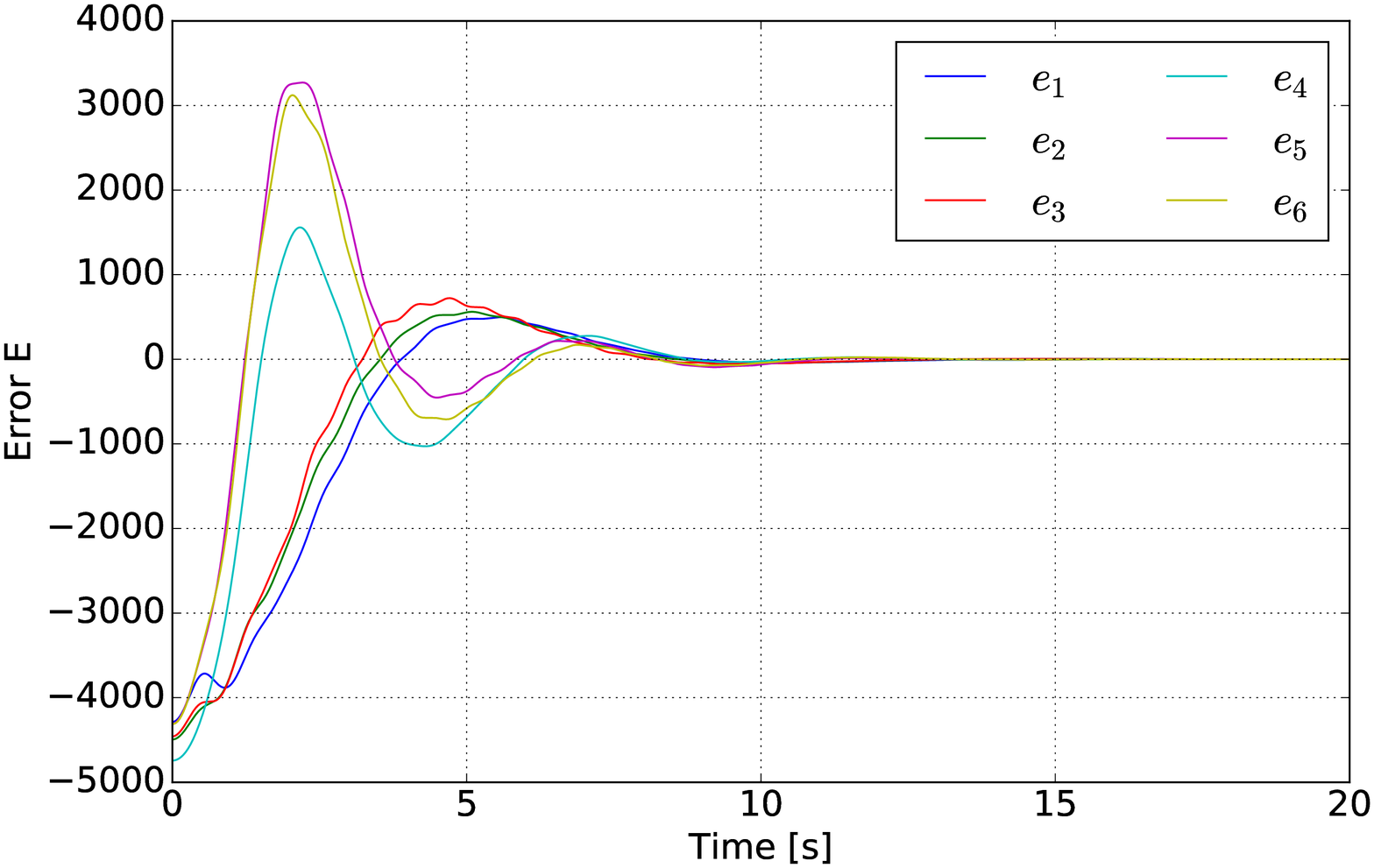}
\caption{}
\end{subfigure}
\begin{subfigure}{0.48\columnwidth}
\includegraphics[width=1\columnwidth]{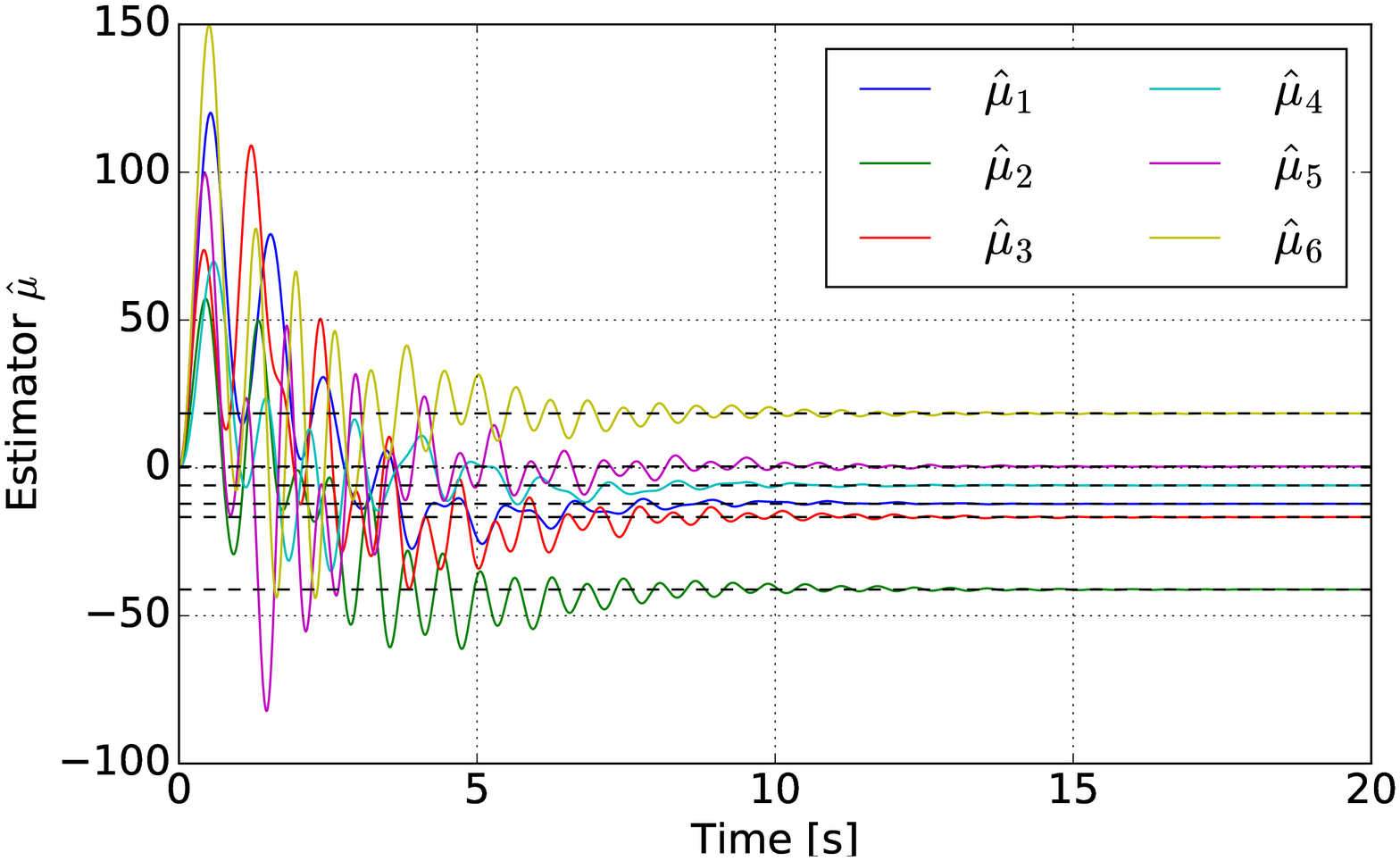}
\caption{}
\end{subfigure}
	\caption{Numerical simulation of a team of four agents with mismatches in their prescribed distances. We employ the estimator proposed in Theorem \ref{th: est}. The final positions of the agents are marked with dots respectively and the red, green, blue and black colors correspond to the agents $1, 2, 3$ and $4$ respectively. The effectiveness of the estimator is shown in (b) and (c), where all the agents' speed and errors (as defined in (\ref{eq: error2})) converge effectively to zero. The plot (d) shows how the six elements of the estimator state $\hat\mu$ converge to the actual $\mu$ (shown in black dashed lines).}
\label{fig: esti}
\end{figure}

We first start validating Theorem \ref{th: est} for a regular tetrahedron formation, with side length equal to $70$ units, whose associated incidence matrix is given by
\begin{equation}
B = \begin{bmatrix}
	-1&  0& -1&  1&  0 & 0\\
	 1&  1&  0&  0&  0 & 1\\
	 0& -1&  1&  0&  1 & 0\\
	 0&  0&  0& -1& -1 & -1
\end{bmatrix}.
\label{eq: Bex1}
\end{equation}
We then randomly generate the following vector of mismatches $\mu$ for each edge $\mathcal{E}_k$
\begin{equation}
	\mu = \begin{bmatrix}12.14 & -41.12 & -16.64 & -5.91 & 0.45 & 18.41\end{bmatrix}^T.
	\label{eq: muex}
\end{equation}
We randomly spread the four agents within an area of $100$ cubic units and with random initial velocities but with speeds smaller than $2$ units per second. We apply the control law as in (\ref{eq: vest}) with the estimator dynamics (\ref{eq: uhat}). We remark that with this setup, one can choose arbitrarily the estimating agents, i.e. how one chooses $B$ for defining the (mismatched) tetrahedron formation does not matter. The results are shown in Figure \ref{fig: esti}.

\begin{figure}
\centering
\begin{subfigure}{0.48\columnwidth}
\includegraphics[width=1\columnwidth]{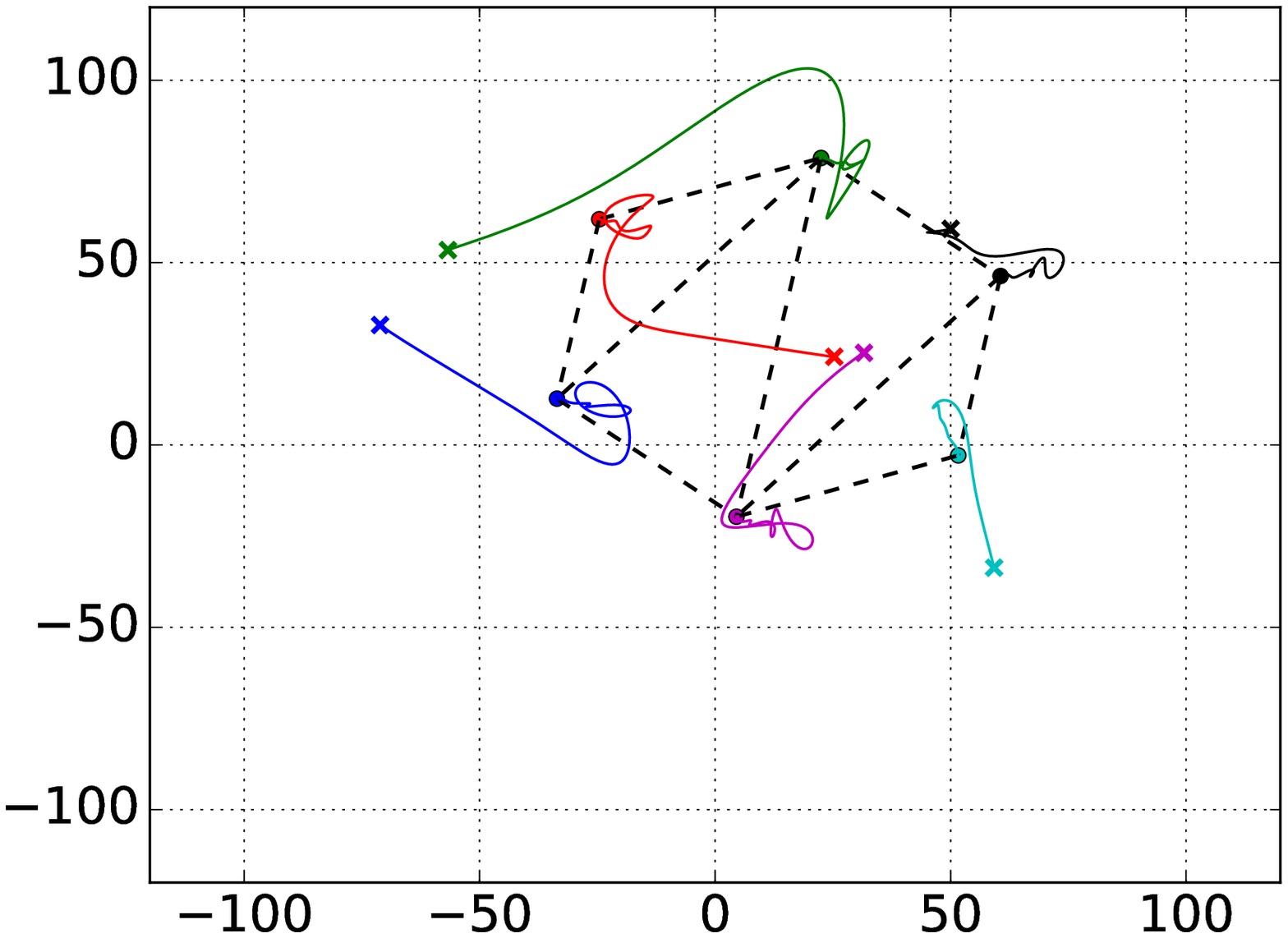}
\caption{}
\end{subfigure}
\begin{subfigure}{0.48\columnwidth}
\includegraphics[width=1\columnwidth]{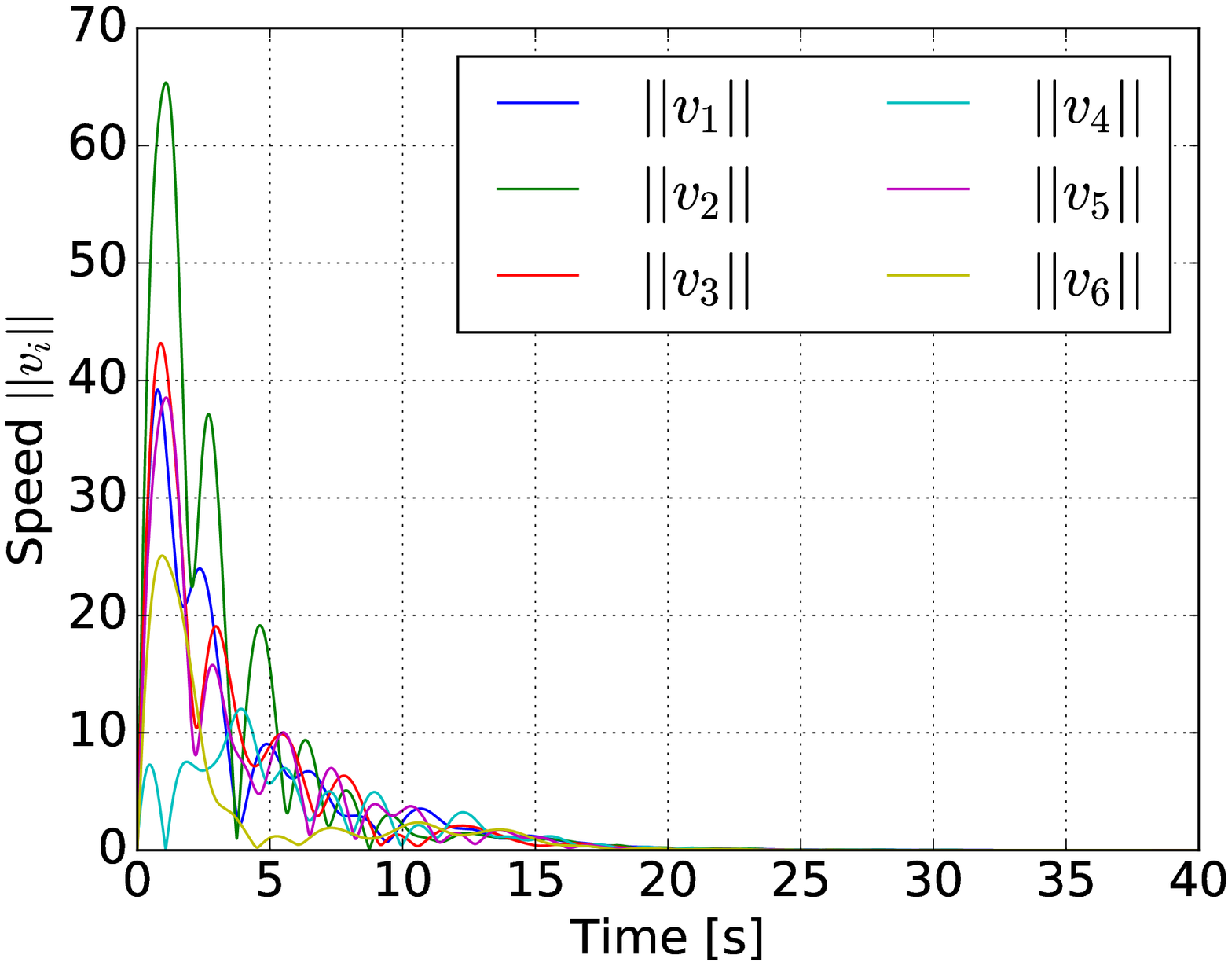}
\caption{}
\end{subfigure}
\begin{subfigure}{0.48\columnwidth}
\includegraphics[width=1\columnwidth]{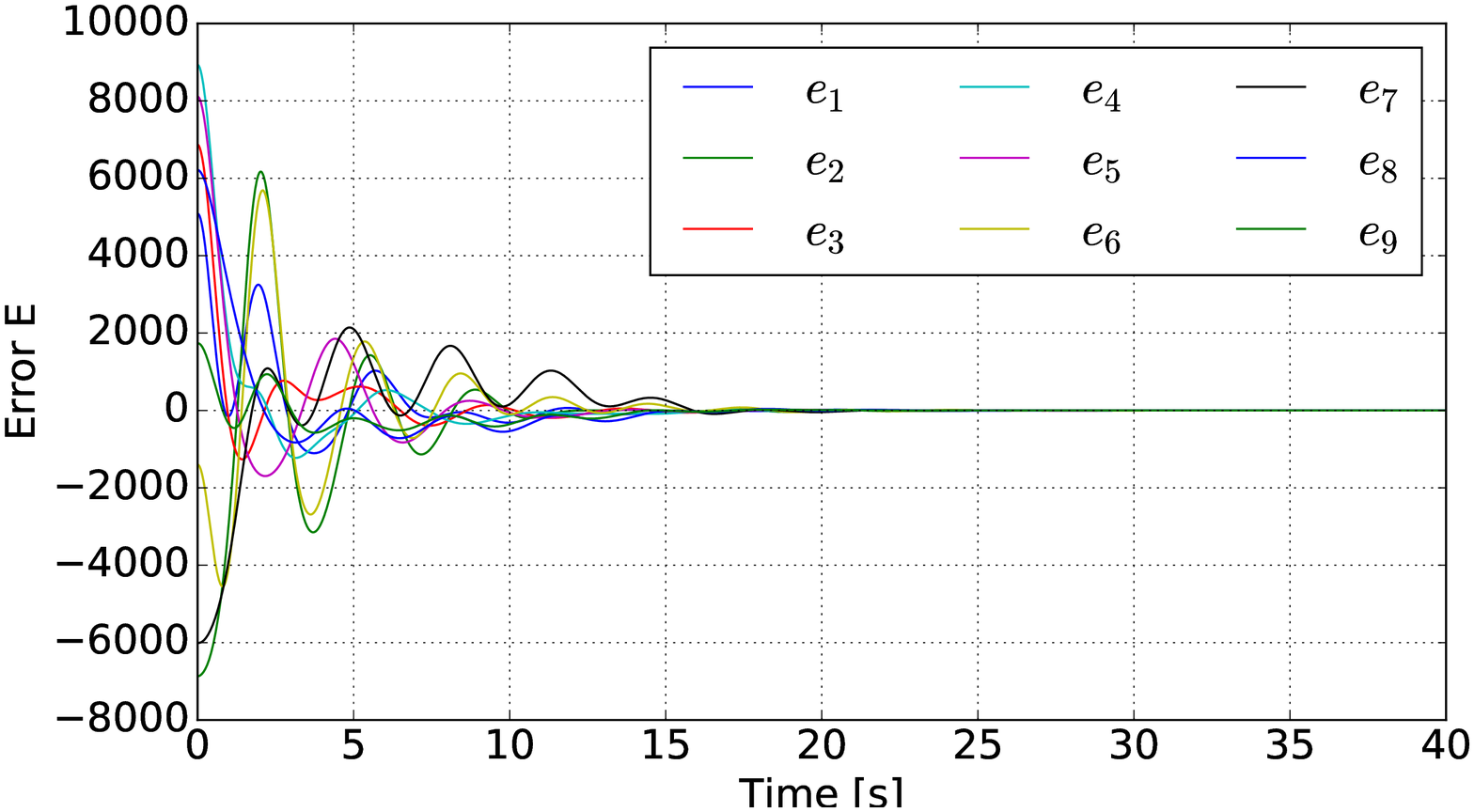}
\caption{}
\end{subfigure}
\begin{subfigure}{0.48\columnwidth}
\includegraphics[width=1\columnwidth]{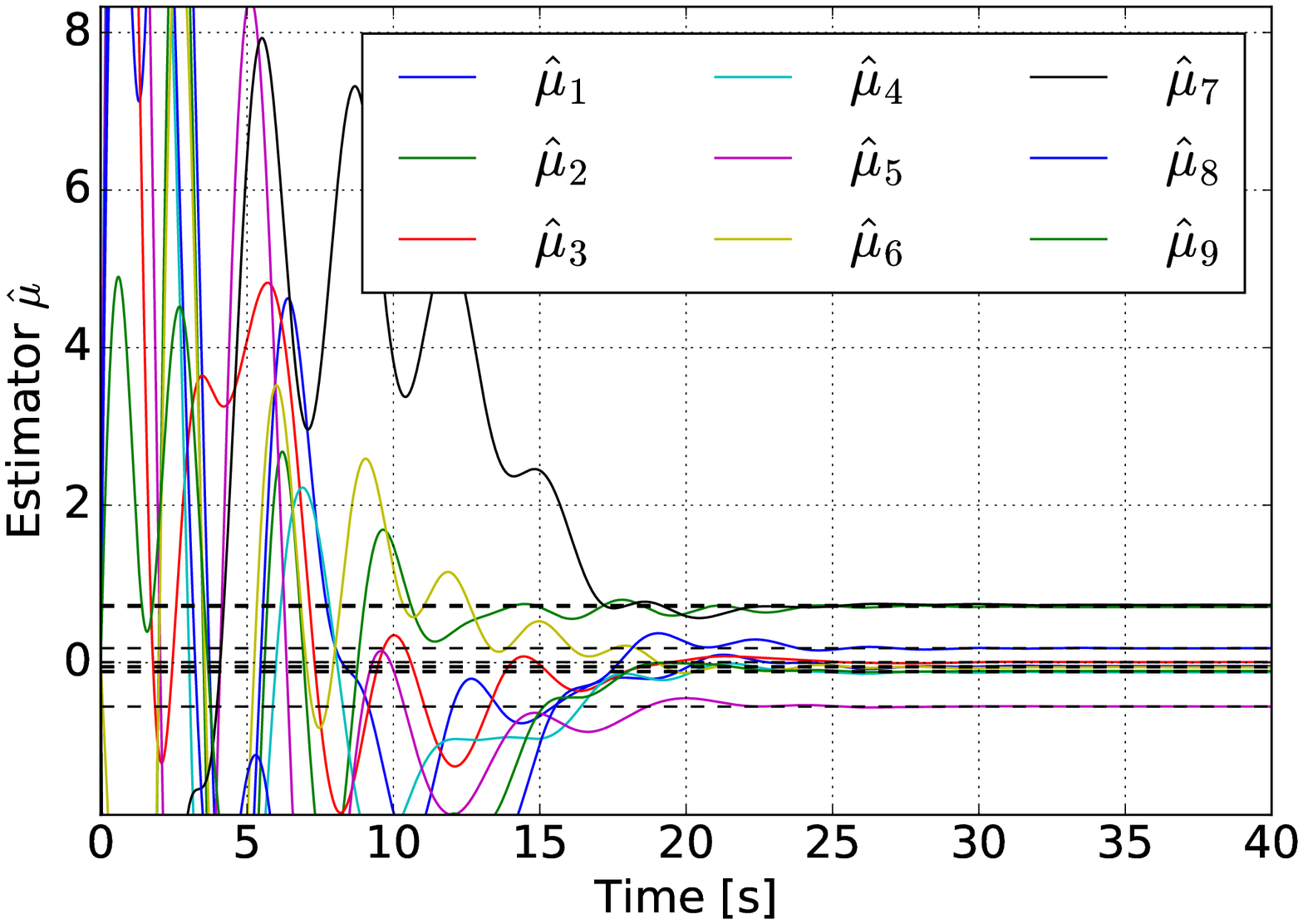}
\caption{}
\end{subfigure}
	\caption{Numerical simulation of a team of six agents for forming a hexagon but with mismatches in their prescribed distances. We employ the estimator proposed in Theorem \ref{th: est2}. The final positions of the agents are marked with dots respectively and the red, green, blue, black, magenta and cyan colors correspond to the agents $1, 2, 3, 4, 5$ and $6$ respectively. The effectiveness of the estimator is shown in (b) and (c), where all the agents' speed and errors (as defined in (\ref{eq: error2})) converge effectively to zero. The plot (d) shows how the nine elements of the estimator state $\hat\mu$ converge to the actual $\mu$ (shown in black dashed lines).}
\label{fig: esti2}
\end{figure}

We now have a team of six agents whose prescribed shape is a regular hexagon, with side length equal to $50$ units, whose incidence matrix is given by
\begin{equation}
B =
\begin{bmatrix}
		1&  0 &  1 &  0 &  0 &  0 &  0 &  0 &  0 \\
       -1&  1&  0&  0&  0&  1&  0&  0&  0 \\
        0& -1& -1&  1&  1&  0&  0&  0&  0 \\
        0&  0&  0& -1&  0&  0&  1&  1&  0 \\
        0&  0&  0&  0& -1& -1& -1&  0&  1 \\
        0&  0&  0&  0&  0&  0&  0& -1& -1
\end{bmatrix},
\end{equation}
and we add the following arbitrary vector of mismatches
\begin{equation}
	\mu = 0.1\begin{bmatrix}-0.43 \\ 7.09 \\ 0.08 \\ -1.19\\ -5.55\\
	-0.574 \\  7.33 \\  1.85 \\ -1.05\end{bmatrix}.
	\label{eq: muhex}
\end{equation}
We use the results in Theorem \ref{th: est2} in order to eliminate the undesired steady-state motion and distortion in the desired shape. We first notice that the estimating agents defined by $S_1$ derived from $B$ are not defining any cycles as it can be checked in Figure \ref{fig: estiHex}. Therefore, the topology for the estimating agents makes Assumption \ref{as: 1} to be satisfied. Since the equilibrium set for the shape is given by $\mathcal{D}$ and not by $\mathcal{Z}$, in order to have a hexagon as a final shape, one restricts the initial positions of the agents to be close to the desired shape and with small initial velocities. We consider the gain $\kappa = 1$ for (\ref{eq: uk}) noting that this value is smaller than the conservative one that can be derived from Theorem \ref{th: est2}. The results are shown in Figure \ref{fig: esti2}.

{\color{black}
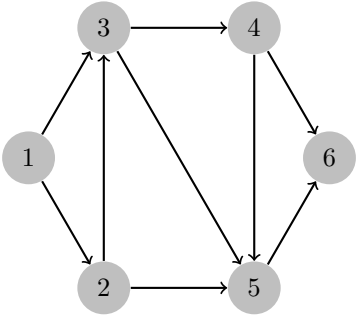
\begin{figure}
\centering
\begin{tikzpicture}[scale=1, auto, swap]
\foreach \pos/\name in {{(-1,-1.73,0)/2}, {(1,-1.73,0)/5}, {(-1,1.73,0)/3}, {(1, 1.73, 0)/4}, 
	{(-2, 0, 0)/1}, {(2, 0, 0)/6}}
        \node[vertex] (\name) at \pos {$\name$};
    \foreach \source/ \dest in {1/2, 2/3, 1/3, 3/4, 3/5, 2/5, 4/5, 4/6, 5/6}
    \path[dedge] (\source) -- (\dest);
\end{tikzpicture}
	\caption{Regular hexagon formation where the tails of the arrows indicate the corresponding estimating agents. Note that this configuration does not contain any cycles, and therefore Assumption \ref{as: 1} is satisfied.}
\label{fig: estiHex}
\end{figure}
}

\begin{figure}
\centering
\includegraphics[width=1\columnwidth]{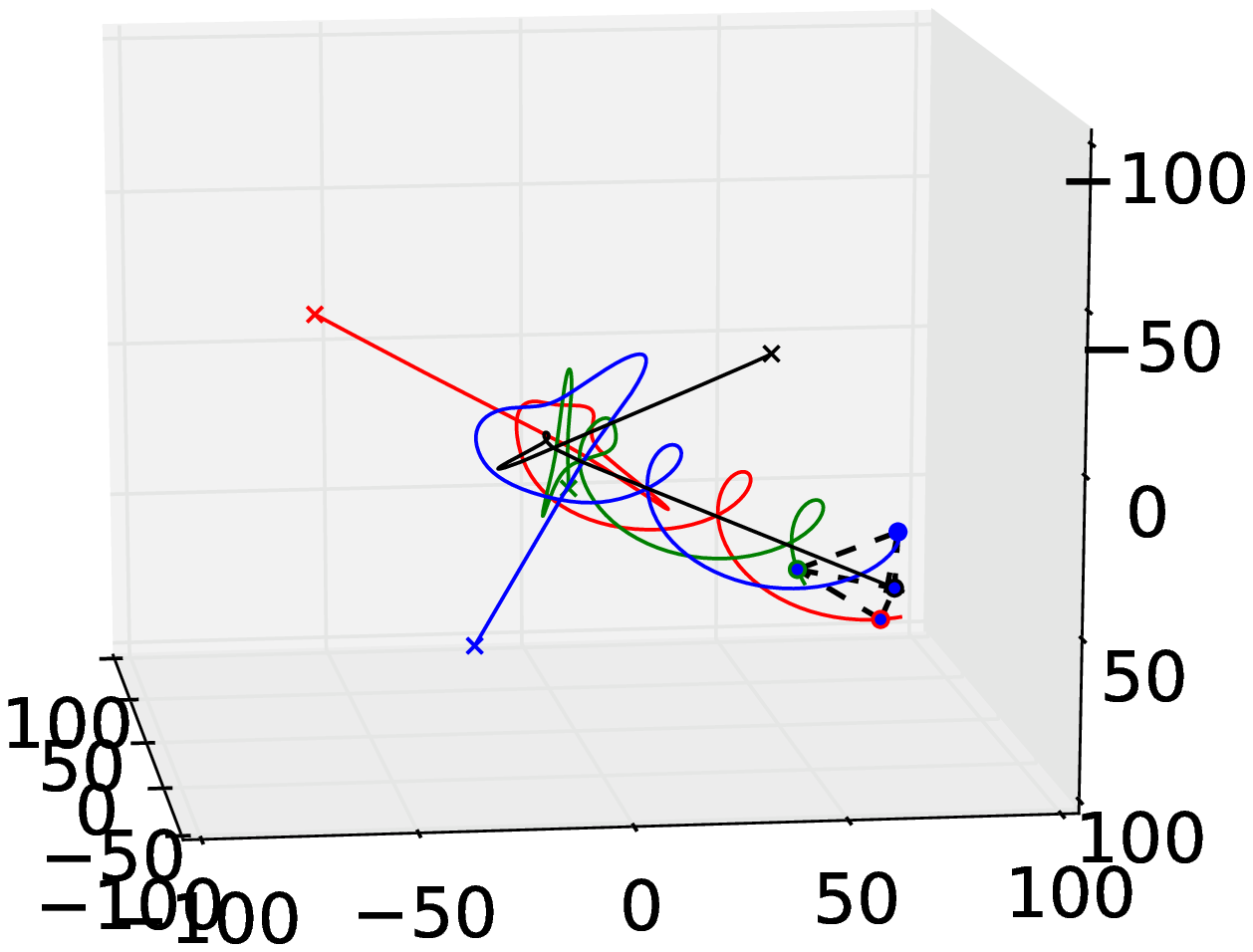}
\begin{subfigure}{0.48\columnwidth}
\includegraphics[width=1\columnwidth]{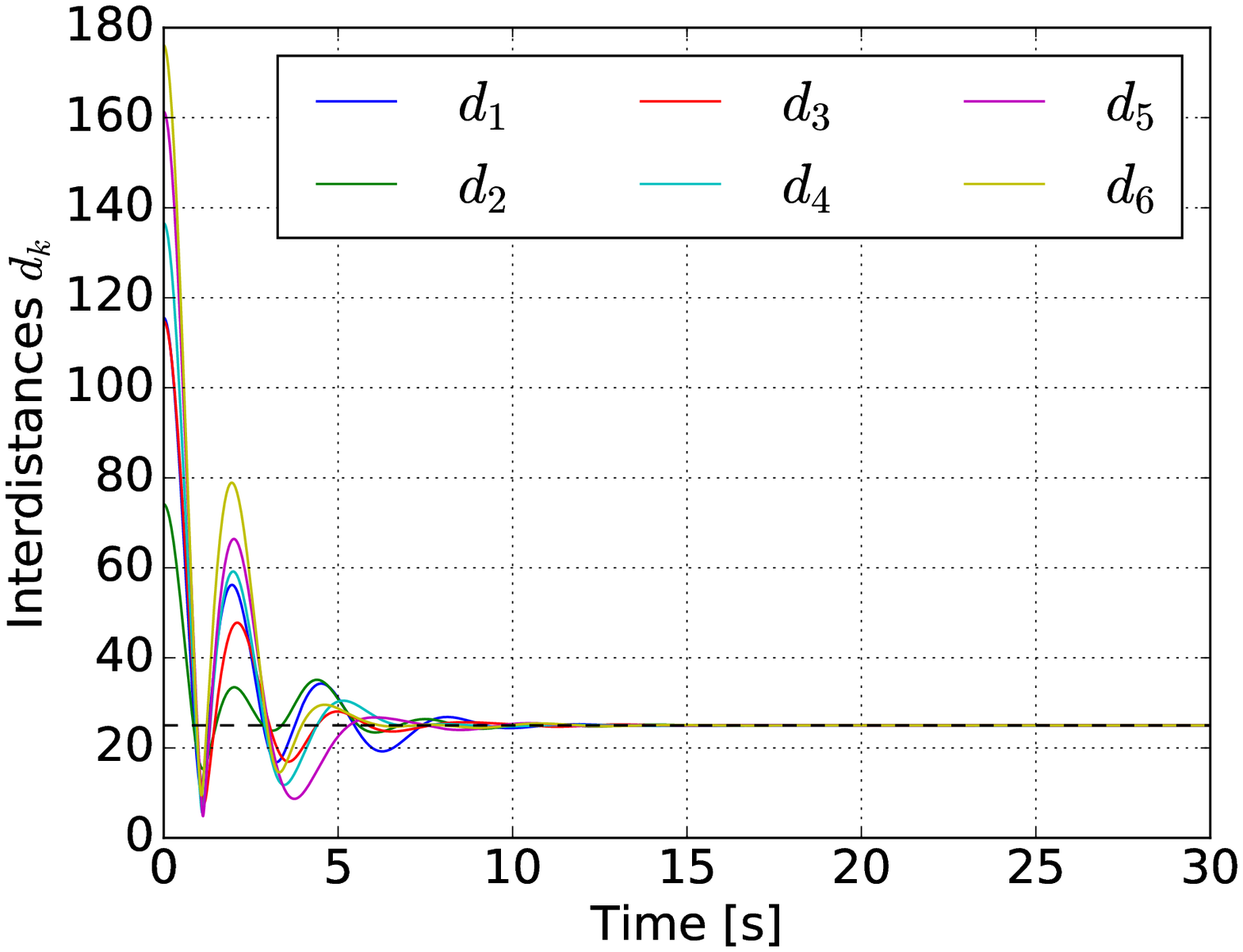}
\caption{}
\end{subfigure}
\begin{subfigure}{0.48\columnwidth}
\includegraphics[width=1\columnwidth]{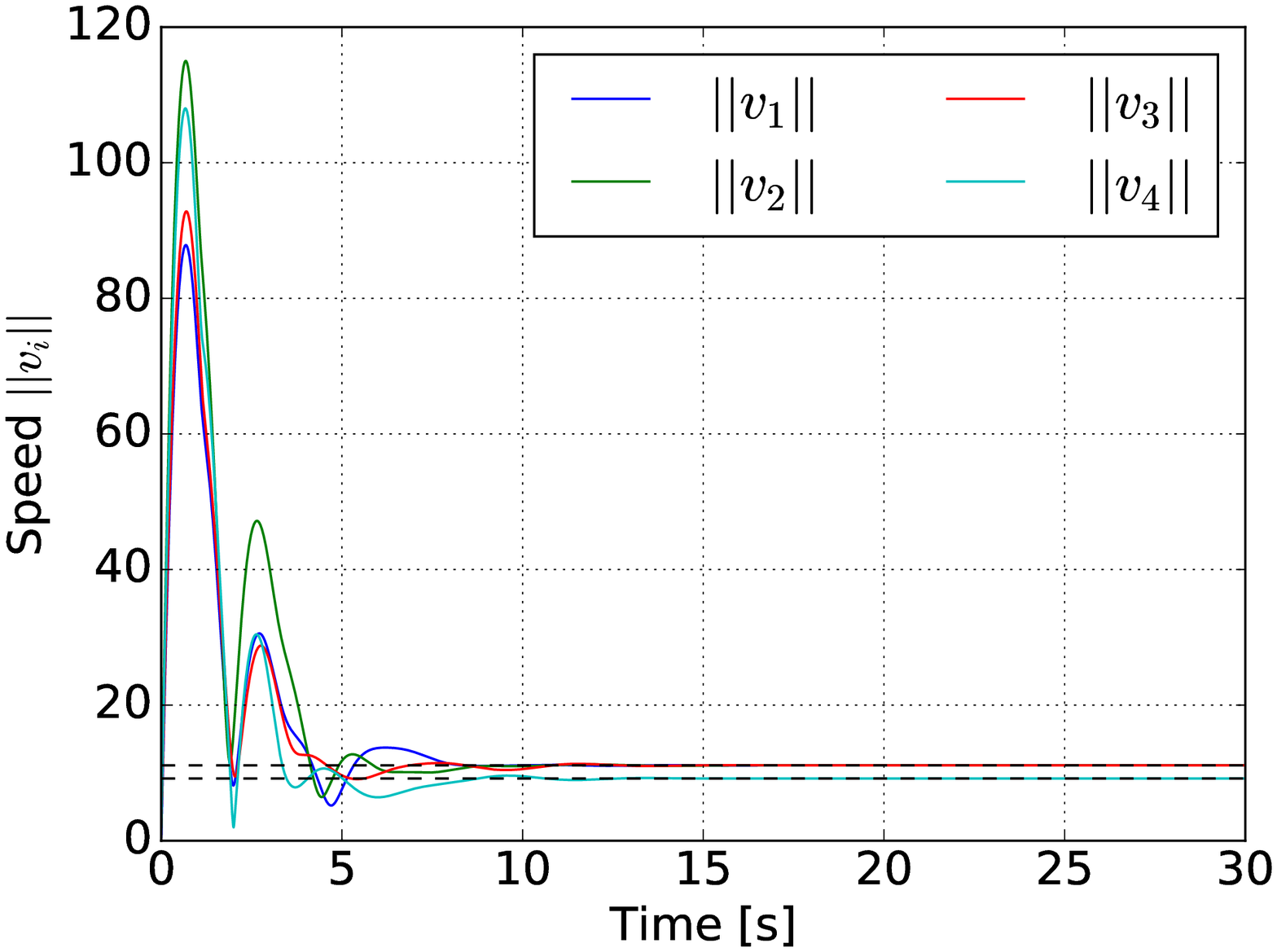}
\caption{}
\end{subfigure}
	\caption{Numerical simulation of a team of four agents travelling in a tetrahedron formation. We employ the results from Theorem \ref{th: rec}. In the first plot the agents $1, 2, 3$ and $4$ are marked with the red, green, blue and black colors respectively. The crosses indicate the initial positions. The black agent \emph{on top} of the tetrahedron follows a linear velocity while the other tree agents follow the same velocity and in addition they are spinning about the centroid of the base of the tetrahedron. We show in plots (a) and (b) the evolution of the inter-agent distances and the speeds of the agents, converging to the desired ones (dashed-lines).}
\label{fig: mot}
\end{figure}

In order to validate the results of Theorem \ref{th: rec} we consider the previous four agents from the first experiment but with a regular tetrahedron of side length equal to $25$ units as a desired shape. The chosen incidence matrix is the same as (\ref{eq: Bex1}). The desired motion of the tetrahedron is given in Figure \ref{fig: motion} but with $^bv_c^*$ being parallel to $^b\omega^*$, i.e., the agent on the \emph{top} of the tetrahedron is following a linear velocity perpendicular to the base, the other three agents follow the same linear velocity, and in addition these three agents also make spinning about the centroid of the base. In order to have such a motion, the distributed motion parameters for $A_v$ as in (\ref{eq: Av}) are given by
\begin{equation}
	\begin{cases}
	\mu_v &= s_v\begin{bmatrix}1 & 1 & 1 & -3 & -3 & -3\end{bmatrix}^T \\
		\tilde\mu_v &= s_v\begin{bmatrix}-1 & -1 & -1 & -1 & -1 & -1\end{bmatrix}^T  \\
\mu_\omega &= s_\omega\begin{bmatrix}1 & 1 & 1 & 0 & 0 & 0\end{bmatrix}^T \\
	\tilde\mu_\omega &= s_\omega\begin{bmatrix}1 & 1 & 1 & 0 & 0 & 0\end{bmatrix}^T, \end{cases}
		\label{eq: mumot}
\end{equation}
where we set $s_v = 0.15$ and $s_\omega = 0.25$, i.e., we regulate the speeds of $^bv_c^*$ and $^b\omega^*$ such that the speed of the agent $4$ is $0.15||z_4^*+z_5^*+z_6^*|| = 9.184$ units per second. Similar calculations can be done for the other three agents in order to derive that their stationary speed will be $11.113$ units per second. We randomly spread the four agents within an area of $50$ cubic units and with random initial velocities with the initial speeds smaller than $2$ units per second. We apply the control law (\ref{eq: uA}) to system (\ref{eq: pdyn}), constructing $A$ with (\ref{eq: mumot}) as in (\ref{eq: A}) and with control gains $c_1 = 1$ and $c_2 = 1$, which are smaller than the conservative ones derived from Theorem \ref{th: rec}, showing the conservative nature of the result of the theorem. The numerical results are shown in Figure \ref{fig: mot}.

%

\section{Conclusions}
\label{sec: con}
	In this paper we have analyzed the effects of having a distance-based controller for rigid formations and second-order agents with the presence of mismatches in their desired inter-agent distances. These effects are a stationary distorted shape with respect to the desired one and an undesired collective motion of the formation. It turns out that both first-order and second-order agents share precisely the same behaviour for the undesired collective motion. We have extended the estimator based solution proposed in \cite{MarCaoJa15} to remove the effects of the mismatches to second-order agents. We have also proposed another estimator with fewer requirements although it only eliminates all the undesired effects for both triangles and tetrahedrons. For the rest of shapes, the new estimator is only effective for removing the undesired state-state collective motion. Nevertheless, a bound on the distortion of the steady-state shape with respect to the desired is given. We have further extended the results from \cite{MaJaCa15} of employing distributed motion parameters in order to control the motion of a desired rigid shape to the second-order agents case. Consequently, it opens possibilities to apply this method directly to actual systems governed by Newtonian dynamics such as quadrotors or marine vessels as it has been shown in \cite{thesis}. We are currently working on extending the recent results for \emph{flexible} formations with mismatches \cite{hecIFAC17a,hecIFAC17b} to second-order agents as well.


%





\ifCLASSOPTIONcaptionsoff
  \newpage
\fi


\bibliographystyle{IEEEtran}
\bibliography{hector_ref}

%




\begin{IEEEbiography}[{\includegraphics[width=1in,height=1.25in,clip,keepaspectratio]{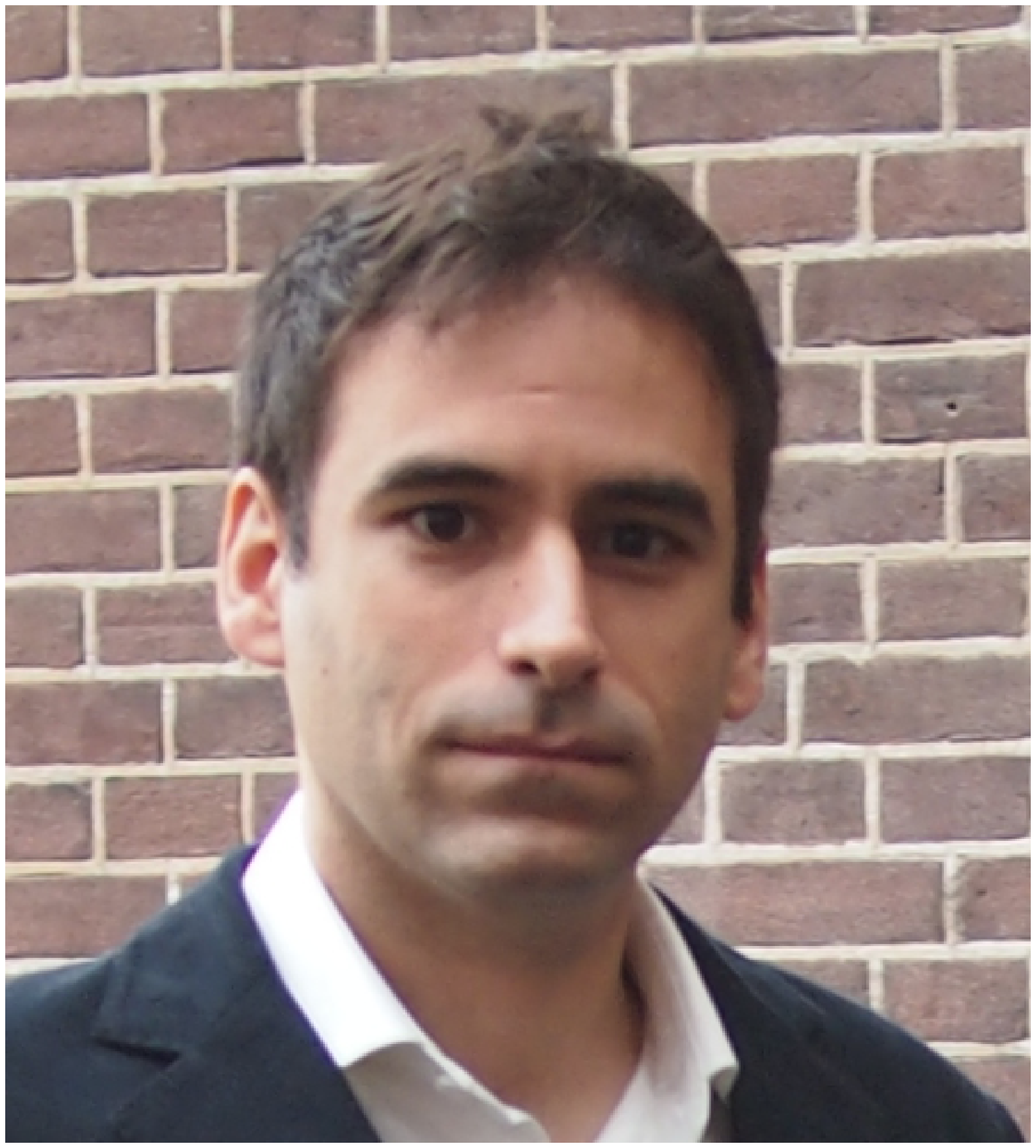}}]{Hector G. de Marina}
	received the M.Sc. degree in electronics engineering from Complutense University of Madrid, Madrid, Spain in 2008 and the M.Sc. degree in control engineering from the University of Alcala de Henares, Alcala de Henares, Spain in 2011. He is a postdoctoral research associate with the Ecole Nationale de l'Aviation Civile (ENAC) in Toulouse, France. His research interests include formation control and navigation for autonomous robots, and drones in particular.
\end{IEEEbiography}
\begin{IEEEbiography}[{\includegraphics[width=1in,height=1.25in,clip,keepaspectratio]{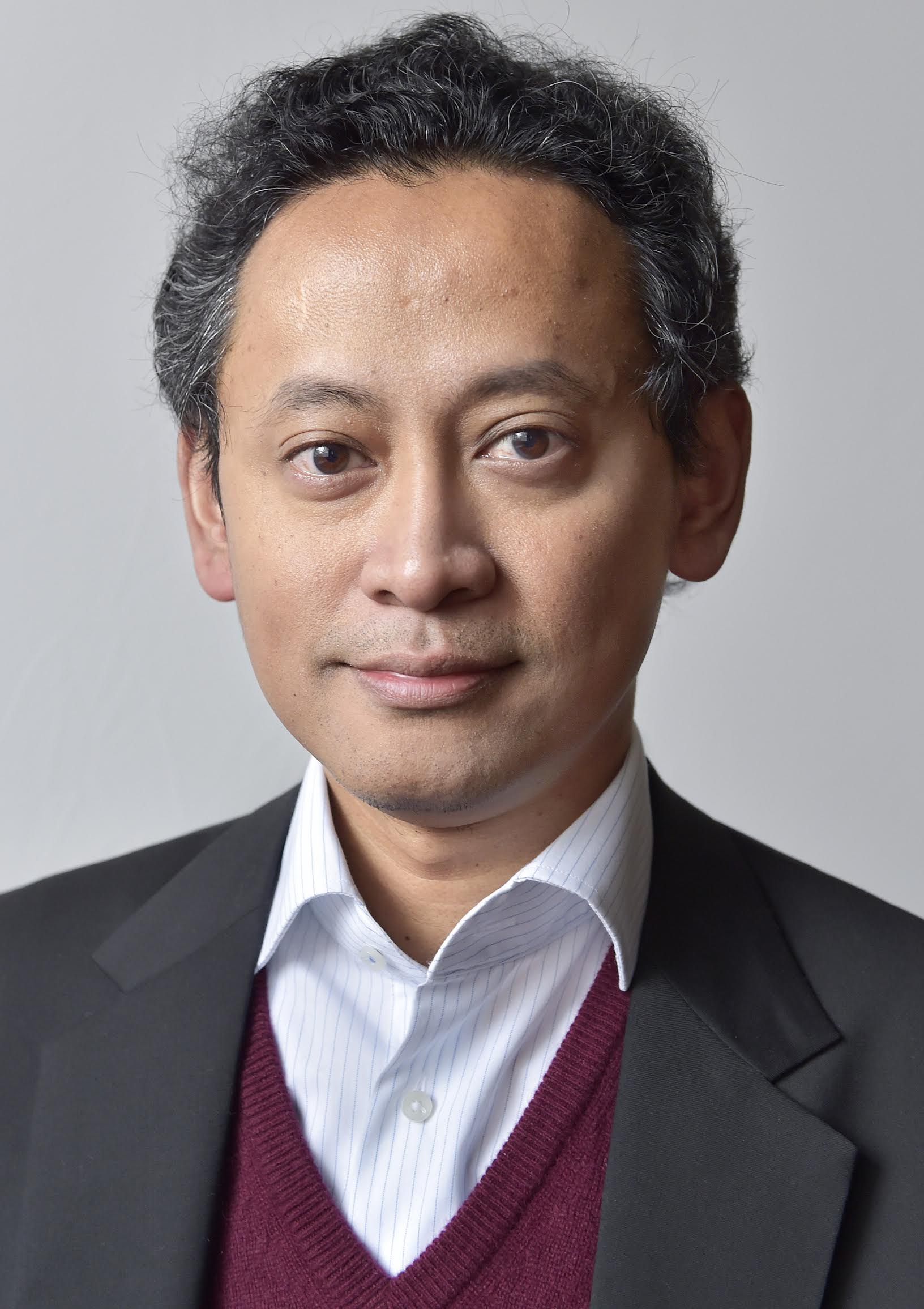}}]{Bayu Jayawardhana}
(SM’13) received the B.Sc. degree in electrical and electronics engineering from the Institut Teknologi Bandung, Bandung, Indonesia, in 2000, the M.Eng. degree in electrical and electronics engineering from the Nanyang Technological University, Singapore, in 2003, and the Ph.D. degree in electrical and electronics engineering from Imperial College London, London, U.K., in 2006. Currently, he is an associate professor in the Faculty of Mathematics and Natural Sciences, University of Groningen, Groningen, The Netherlands. He was with Bath University, Bath, U.K., and with Manchester Interdisciplinary Biocentre, University of Manchester, Manchester, U.K. His research interests are on the analysis of nonlinear systems, systems with hysteresis, mechatronics, systems and synthetic biology. 
Prof. Jayawardhana is a Subject Editor of the International Journal of Robust and Nonlinear Control, an associate editor of the European Journal of Control, and a member of the Conference Editorial Board of the IEEE Control Systems Society.
\end{IEEEbiography}
\begin{IEEEbiography}[{\includegraphics[width=1in,height=1.25in,clip,keepaspectratio]{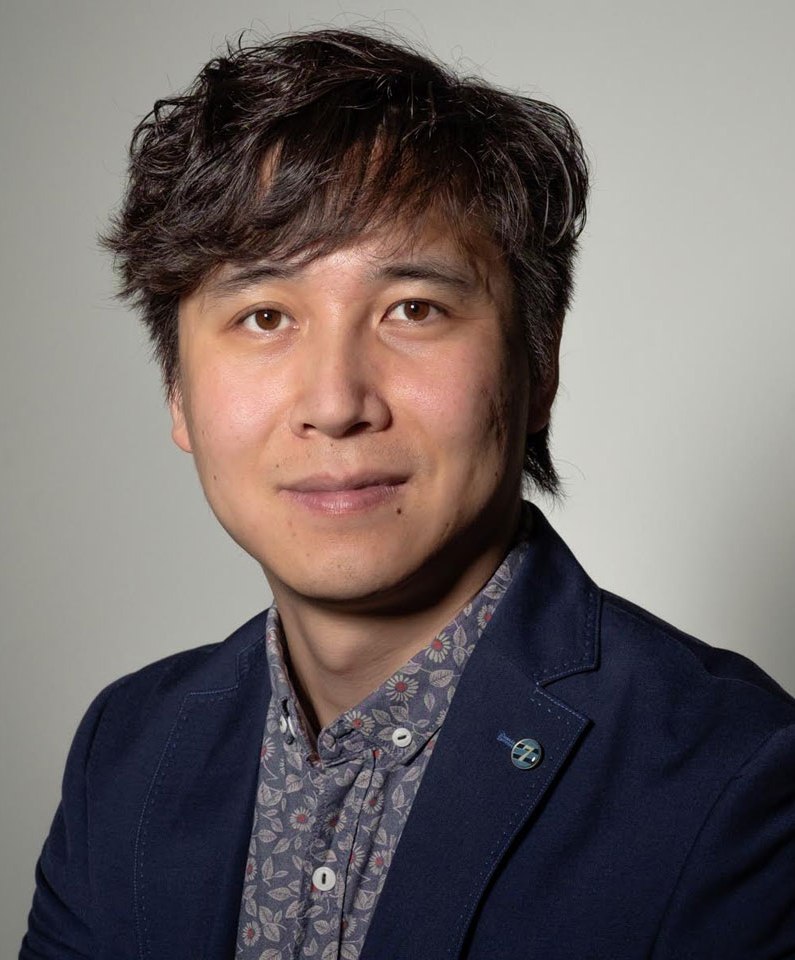}}]{Ming Cao}
is currently professor of systems and control with the Engineering and Technology Institute (ENTEG) at the University of Groningen, the Netherlands, where he started as a tenure-track assistant professor in 2008. He received the Bachelor degree in 1999 and the Master degree in 2002 from Tsinghua University, Beijing, China, and the PhD degree in 2007 from Yale University, New Haven, CT, USA, all in electrical engineering. From September 2007 to August 2008, he was a postdoctoral research associate with the Department of Mechanical and Aerospace Engineering at Princeton University, Princeton, NJ, USA. He worked as a research intern during the summer of 2006 with the Mathematical Sciences Department at the IBM T. J. Watson Research Center, NY, USA.  He is the 2017 recipient of the Manfred Thoma medal from the International Federation of Automatic Control (IFAC) and the 2016 recipient of the European Control Award from the European Control Association (EUCA). He is an associate editor for IEEE Transactions on Automatic Control, IEEE Transactions on Circuits and Systems and Systems and Control Letters, and for the Conference Editorial Board of the IEEE Control Systems Society. He is also a member of the IFAC Technical Committee on Networked Systems. His main research interest is in autonomous agents and multi-agent systems, mobile sensor networks and complex networks.
\end{IEEEbiography}

\end{document}

%% file: FIGsquaredef.tex
\begin{tikzpicture}[line join=round]
\draw(0,0)--(1,0)--(1,1)--(0,1)--(0,0);
\filldraw(0,0) circle (2pt);
\filldraw(1,0) circle (2pt);
\filldraw(1,1) circle (2pt);
\filldraw(0,1) circle (2pt);
\filldraw(-.85,.85) circle (2pt);
\filldraw(.15,.85) circle (2pt);
\draw(0,0)--(-.85,.85);
\draw(1,0)--(.15,.85);
\draw(-.85,.85)--(.15,.85);
\draw[red,->,style=dashed] (-.85,.85) to [bend left=5] (0,1);\draw[red,->,style=dashed] (.15,.85) to [bend left=5] (1,1);\end{tikzpicture}

%% file: FIGsquareRig.tex
\begin{tikzpicture}[line join=round]
\draw(0,0)--(1,0)--(1,1)--(0,1)--(0,0)--(1,1);
\filldraw(0,0) circle (2pt);
\filldraw(1,0) circle (2pt);
\filldraw(1,1) circle (2pt);
\filldraw(0,1) circle (2pt);
\end{tikzpicture}

%% file: FIGtetra.tex
\begin{tikzpicture}[line join=round]
[\tikzset{>=latex}]\filldraw[draw=black,fill=white,fill opacity=0.8](-.408,.224)--(0,.625)--(0,.625)--(.68,.032)--cycle;
\filldraw[draw=black,fill=white,fill opacity=0.8](-.273,-.256)--(0,.625)--(0,.625)--(-.408,.224)--cycle;
\filldraw[draw=black,fill=white,fill opacity=0.8](.68,.032)--(0,.625)--(0,.625)--(-.273,-.256)--cycle;
\end{tikzpicture}

%% file: FIGmotionTetra.tex
\begin{tikzpicture}[line join=round]
[\tikzset{>=latex}]\filldraw[fill=white](-2.353,1.014)--(-1.808,-.927)--(2.955,-.705)--(2.41,1.236)--cycle;
\filldraw[draw=black,fill=white,fill opacity=0.8](.55,.197)--(.958,.598)--(.958,.598)--(1.638,.005)--cycle;
\filldraw[draw=black,fill=white,fill opacity=0.8](.685,-.283)--(.958,.598)--(.958,.598)--(.55,.197)--cycle;
\draw[draw=red,arrows=->](2.381,.111)--(2.381,1.112);
\draw[draw=black,arrows=<->,thick](-.368,-.243)--(-.942,.059)--(-.569,.51);
\draw[arrows=-,thick](.992,.009)--(.958,.129)--(1.468,.153);
\draw[draw=black,arrows=-](.958,.129)--(.958,.598);
\draw[draw=black,arrows=-](.958,.129)--(1.148,.279);
\draw[arrows=-,thick](.958,.129)--(.958,.598);
\draw[draw=black,arrows=-](.958,.129)--(1.309,.01);
\filldraw[draw=black,fill=white,fill opacity=0.8](1.638,.005)--(.958,.598)--(.958,.598)--(.685,-.283)--cycle;
\draw[arrows=->,thick](1.468,.153)--(2.115,.183);
\draw[draw=black,arrows=->](1.309,.01)--(1.847,-.172);
\draw[draw=black,arrows=-](1.148,.279)--(1.378,.46);
\draw[thick](1.162,1.08)--(1.167,1.082)--(1.171,1.083)--(1.176,1.085)--(1.18,1.086)--(1.185,1.088)--(1.189,1.09)--(1.194,1.091)--(1.198,1.093)--(1.202,1.095)--(1.206,1.097)--(1.211,1.098)--(1.215,1.1)--(1.218,1.102)--(1.222,1.104)--(1.226,1.106)--(1.23,1.108)--(1.233,1.11)--(1.237,1.112)--(1.24,1.114)--(1.244,1.116)--(1.247,1.118)--(1.25,1.12)--(1.253,1.122)--(1.256,1.124)--(1.259,1.126)--(1.262,1.128)--(1.264,1.13)--(1.267,1.133)--(1.27,1.135)--(1.272,1.137)--(1.274,1.139)--(1.277,1.141)--(1.279,1.144)--(1.281,1.146)--(1.283,1.148)--(1.285,1.151)--(1.286,1.153)--(1.288,1.155)--(1.29,1.158)--(1.291,1.16)--(1.292,1.162)--(1.294,1.165)--(1.295,1.167)--(1.296,1.17)--(1.297,1.172)--(1.297,1.174)--(1.298,1.177)--(1.299,1.179)--(1.299,1.182)--(1.3,1.184)--(1.3,1.187)--(1.3,1.189)--(1.3,1.191)--(1.3,1.194)--(1.3,1.196)--(1.3,1.199)--(1.3,1.201)--(1.299,1.204)--(1.299,1.206)--(1.298,1.208)--(1.297,1.211)--(1.297,1.213)--(1.296,1.216)--(1.295,1.218)--(1.293,1.22)--(1.292,1.223)--(1.291,1.225)--(1.289,1.228)--(1.288,1.23)--(1.286,1.232)--(1.284,1.235)--(1.283,1.237)--(1.281,1.239)--(1.279,1.242)--(1.277,1.244)--(1.274,1.246)--(1.272,1.248)--(1.269,1.251)--(1.267,1.253)--(1.264,1.255)--(1.262,1.257)--(1.259,1.259)--(1.256,1.261)--(1.253,1.263)--(1.25,1.266)--(1.247,1.268)--(1.243,1.27)--(1.24,1.272)--(1.237,1.274)--(1.233,1.276)--(1.23,1.278)--(1.226,1.28)--(1.222,1.281)--(1.218,1.283)--(1.214,1.285)--(1.21,1.287)--(1.206,1.289)--(1.202,1.29)--(1.198,1.292)--(1.194,1.294)--(1.189,1.296)--(1.185,1.297)--(1.18,1.299)--(1.176,1.3)--(1.171,1.302)--(1.166,1.303)--(1.162,1.305)--(1.157,1.306)--(1.152,1.308)--(1.147,1.309)--(1.142,1.31)--(1.137,1.312)--(1.132,1.313)--(1.126,1.314)--(1.121,1.315)--(1.116,1.316)--(1.111,1.318)--(1.105,1.319)--(1.1,1.32)--(1.094,1.321)--(1.089,1.322)--(1.083,1.322)--(1.078,1.323)--(1.072,1.324)--(1.066,1.325)--(1.061,1.326)--(1.055,1.326)--(1.049,1.327)--(1.044,1.328)--(1.038,1.328)--(1.032,1.329)--(1.026,1.329)--(1.02,1.33)--(1.014,1.33)--(1.008,1.331)--(1.002,1.331)--(.997,1.331)--(.991,1.332)--(.985,1.332)--(.979,1.332)--(.973,1.332)--(.967,1.332)--(.961,1.332)--(.955,1.332)--(.949,1.332)--(.943,1.332)--(.937,1.332)--(.931,1.332)--(.925,1.332)--(.919,1.331)--(.913,1.331)--(.907,1.331)--(.901,1.33)--(.895,1.33)--(.89,1.329)--(.884,1.329)--(.878,1.328)--(.872,1.328)--(.866,1.327)--(.861,1.326)--(.855,1.326)--(.849,1.325)--(.843,1.324)--(.838,1.323)--(.832,1.322)--(.827,1.322)--(.821,1.321)--(.816,1.32)--(.81,1.319)--(.805,1.317)--(.8,1.316)--(.794,1.315)--(.789,1.314)--(.784,1.313)--(.779,1.312)--(.774,1.31)--(.769,1.309)--(.764,1.308)--(.759,1.306)--(.754,1.305)--(.749,1.303)--(.745,1.302)--(.74,1.3)--(.735,1.299)--(.731,1.297)--(.726,1.295)--(.722,1.294)--(.718,1.292)--(.714,1.29)--(.709,1.289)--(.705,1.287)--(.701,1.285)--(.697,1.283)--(.694,1.281)--(.69,1.279)--(.686,1.277)--(.683,1.276)--(.679,1.274)--(.676,1.272)--(.672,1.27)--(.669,1.268)--(.666,1.265)--(.663,1.263)--(.66,1.261)--(.657,1.259)--(.654,1.257)--(.651,1.255)--(.649,1.253)--(.646,1.25)--(.644,1.248)--(.642,1.246)--(.639,1.244)--(.637,1.241)--(.635,1.239)--(.633,1.237);
\draw[draw=black,dashed](1.847,.829)--(1.847,-.172);
\draw[arrows=->,thick](-.942,.059)--(-.932,.476);
\draw[draw=black,arrows=->](.958,.598)--(.958,.942);
\draw[arrows=->,thick](.958,.598)--(.958,1.193);
\draw[arrows=<-,thick](1.09,-.342)--(.992,.009);
\draw[draw=black,arrows=->](1.378,.46)--(1.847,.829);
\draw[draw=black,dashed](1.847,.829)--(.958,.942);
\draw[thick,->](1.206,1.097)--(1.055,1.059);
\node at (-1.233,.288) {$O_g$};\node at (.428,.484) {$O_b$};\node at (.958,1.568) {$^b\omega^*$};\node at (2.012,.878) {$^bv^*_c$};\draw[red,->] (.958,.129) to [bend left=45] (1.631,-.897);\end{tikzpicture}